\documentclass[11pt]{amsart}

\usepackage{a4wide}
\usepackage{bbm}
\usepackage{t1enc} 
\usepackage{epsf}
\usepackage{amssymb}
\usepackage{graphicx}
\usepackage{psfrag}

\newtheorem{theorem}{Theorem}
\newtheorem{lemma}{Lemma}
\newtheorem{proposition}{Proposition}
\newtheorem{corollary}{Corollary}

\theoremstyle{definition}
\newtheorem{definition}{Definition}
\newtheorem{remark}{Remark}

\begin{document}

\title[Discrete Tomography]{Discrete Tomography of Planar Model Sets}

\author[Baake, Gritzmann, Huck, Langfeld, and Lord]{M. Baake}

\author[]{P. Gritzmann}
\author[]{C. Huck}
\author[]{B. Langfeld}
\author[]{K. Lord}
\address{\hspace*{-1em} Fakult\"{a}t f\"{u}r Mathematik,
  Universit\"{a}t Bielefeld, Postfach 10 01 31, 33501 Bielefeld, Germany}
\email{\{mbaake,huck\}@math.uni-bielefeld.de}
\urladdr{http://www.math.uni-bielefeld.de/baake/}

\address{\hspace*{-1em} Zentrum Mathematik, TU M\"{u}nchen,
 Boltzmannstr. 3, 85747 Garching bei M\"{u}nchen, Germany}
\email{\{gritzman,langfeld,lord\}@ma.tum.de}
\urladdr{http://www-m9.ma.tum.de/dm/homepages/gritzmann/}

\begin{abstract} 
  Discrete tomography is a well-established method to investigate
  finite point sets, in particular finite subsets of periodic systems.
  Here, we start to develop an efficient approach for the treatment of
  finite subsets of mathematical quasicrystals. To this end, the class
  of cyclotomic model sets is introduced, and the corresponding
  consistency, reconstruction and uniqueness problems of the discrete
  tomography of these sets are discussed.
\end{abstract}

\maketitle

\section{Introduction}
\label{intro}

\emph{Discrete tomography} is concerned with the inverse problem of
retrieving information about some discrete object from (generally
noisy) information about its incidences with certain query sets.  A
typical example is the \emph{reconstruction} of a finite point set
from its line sums in a small number $m$ of directions. The term
$X$-ray (or $X$-ray projection) is a {\em generic} name here which
stands for a mechanism that produces {\em weighted} projection data.
More precisely, a (discrete parallel) $X$-ray of a finite subset of
Euclidean $d$-space $\mathbbm{R}^d$ in direction $u$ gives the number
of points in the set on each line in $\mathbbm{R}^d$ parallel to $u$.
(This concept should not be misunderstood in the sense of diffraction
theory, where $X$-rays provide rather different information on the
underlying structure that is based on statistical pair correlations;
compare with~Guinier (1994), Cowley (1995) and Fewster (2003).)

Many papers focus on the discrete tomography of subsets of lattices
since lattices are good models for crystalline structures.  However,
nature provides us also with \emph{structured} non-lattice sets, the
so-called \emph{quasicrystals}. In the present paper, we shall
investigate the discrete tomography of systems of {\em aperiodic
  order}, more precisely, of so-called \emph{model sets} (or
\emph{mathematical quasicrystals}), which are commonly accepted as a
mathematical model for perfect quasicrystalline structures in
nature~(Steurer, 2004).  As model sets possess a `dimensional
hierarchy', which means that any model set in $d$ dimensions can be
sliced into model sets of dimension $d-1$, solving the reconstruction
problem for two-dimensional systems with aperiodic order lies at the
heart of solving the corresponding problem in three dimensions.

The main motivation for our interest in the discrete tomography of
model sets comes from the demand of materials science to reconstruct
three-dimensional (quasi)crystals or planar layers of them from their
images obtained with quantitative {\em high resolution transmission
  electron microscopy} (HRTEM) in a small number of directions.

In fact, in~Schwander {\it et al.} (1993) and Kisielowski {\it et al.} (1995),
the technique QUANTITEM (quantitative analysis of the information
coming from transmission electron microscopy) is described, which is
based on HRTEM and can effectively measure the number of atoms lying
on lines parallel to certain directions. At present, the measurement
of the number of atoms lying on a line can only be achieved for some
crystals; see~Schwander {\it et al.} (1993) and Kisielowski {\it et
  al.} (1995). However, it is reasonable to expect that future developments in
technology will improve this situation.

Roughly speaking, planar model sets are projections of certain subsets
depending on some \emph{window} $W$ of a higher dimensional lattice
into the plane. In Section \ref{sec3} we will define model sets in
general, but we will mainly restrict ourselves to a well-known class
of \emph{planar} model sets, the \emph{cyclotomic model sets}.  On the
one hand, cyclotomic model sets exhibit a particularly nice and useful
algebraic structure, while on the other hand real-world quasicrystals
can be sliced into parallel planar layers that can be modeled by
cyclotomic model sets~(Pleasants, 2000).  Also, in a certain sense,
cyclotomic model sets can be seen as a direct generalization of the
square lattice $\mathbbm{Z}^2$, the classical planar setting of
discrete tomography.

Naturally, all classic issues of discrete tomography including
uniqueness, reconstruction and stability (see {\it e.g.} the book by Herman \&
Kuba (1999) and, in particular, the papers by 
 Gardner \& Gritzmann (1997), (Gritzmann {\it et al.} (1998), Gardner
 {\it et
al.} (1999), Gritzmann {\it et al.} (2000), Alpers {\it et al.} (2001) and Alpers \& Gritzmann
 (2006)) can be studied for model sets
as, in principle, they are just different ground sets, for the
potential solutions.  As it turns out, however, the more general
setting does disclose some new aspects, and the present paper will
stress these. In particular, it is a priori not even clear how to
decide whether a translate of a given finite point set occurs within
an aperiodic structure.

As a matter of fact, previous studies have focussed on the `anchored'
case that the underlying ground set is located in a linear space,
i.e., in a space with a specified location of the origin.  The $X$-ray
data is then taken with respect to this localization.  This assumption
is mainly justified by the fact that, as point sets, one has the
equality $t+\mathbbm{Z}^2=\mathbbm{Z}^2$ for all $t\in\mathbbm{Z}^2$.
Hence, in the lattice case, one can always assume that -- if a
solution exists -- it is close to the origin. In the affine and
aperiodic case of planar model sets it is a priori not clear how far
out solutions may exist and how one can systematically search for
them.

The main result of this paper is, however, that for cyclotomic model
sets (coming from polyhedral windows) all possible localizations can
be determined efficiently.  In fact, we shall solve a corresponding
{\em decomposition problem} and a {\em separation problem}. This will
allow us to reduce tomographic problems such as reconstruction and
uniqueness for cyclotomic model sets to the corresponding classical
problems with certain restrictions. One difference is manifest in the
fact that potential solutions are subsets of a finite list of patches,
whose number typically grows polynomially in the size. In fact, using
the algebraic and the geometric structure of cyclotomic model sets we
show that in a well-defined way the algorithmic methods that have been
developed for the lattice case can be extended to the discrete
tomography of cyclotomic model sets.  (Note, however, as a warning
that even in the (linear) lattice case $\mathbbm{Z}^2$ these problems
are $\mathbbm{NP}$-hard for three or more lattice directions;
see~Gritzmann {\it et al.} (1998) and Gardner {\it et al.} (1999).)

Let us be more specific.  By using the Minkowski representation of
algebraic number fields, we introduce, for $n \notin \{1,2\}$, the
corresponding class of \emph{cyclotomic model sets}
$\varLambda\subset\mathbbm{C}\cong\mathbbm{R}^2$ which live on
$\mathbbm{Z}[\zeta_{n}]\subset\mathbbm{C}$, where $\zeta_{n}$ is a
primitive $n$th root of unity in $\mathbbm{C}$, {\it e.g.},
$\zeta_{n}=e^{2\pi i/n}$.  (Here, and in the following a subset $S$ of
$\mathbbm{R}^{2}$ is said to \emph{live on} a subgroup $G$ of
$\mathbbm{R}^{2}$ if its difference set $S-S:=\{s-s'\,|\,s,s'\in S\}$
is a subset of $G$. Obviously, this is equivalent to the existence of
a suitable $t\in \mathbbm{R}^{2}$ such that $S\subset t+G$.)  The
$\mathbbm{Z}$-module $\mathbbm{Z}[\zeta_{n}]$ is the ring of integers
in the $n$th cyclotomic field $\mathbbm{Q}(\zeta_{n})$, and, for $n
\notin \{1,2,3,4,6\}$, when viewed as a subset of the plane, is dense;
see Section \ref{sec2} for details.  In contrast, (cyclotomic) model
sets $\varLambda$ are \emph{Delone sets}, i.e., they are uniformly
discrete and relatively dense. In fact, model sets are even {\em
  Meyer sets}, meaning that also $\varLambda-\varLambda$ is uniformly
discrete; see~(Moody, 2000). It turns out that, excepting the
cyclotomic model sets living on $\mathbbm{Z}[\zeta_{n}]$ with $n \in
\{3,4,6\}$ (these are exactly the translations of the square and the
triangular lattice, respectively), cyclotomic model sets $\varLambda$
are \emph{aperiodic}, meaning that they have no translational
symmetries. Well-known examples with $N$-fold cyclic symmetry are the
vertex sets of the square tiling ($n=N=4$), the triangle tiling
($2n=N=6$), the Ammann-Beenker tiling ($n=N=8$), the T\"ubingen
triangle tiling ($2n=N=10$) and the shield tiling ($n=N=12$),
respectively; see below for details. Observe that $5,8,10$ and $12$
are standard cyclic symmetries of genuine planar
quasicrystals~(Steurer, 2004).

Whether or not one has future applications in materials science of
quasicrystals in mind, the starting point will always be a specific
structure model. This means that the specific type of the
(quasi)crystal is known, and one is confronted with the $X$-ray data
of an unknown finite subset of it. Let us point out that the
rotational orientation of the probe in an electron microscope can
rather easily be ascertained in the diffraction mode, prior to taking
images in the high-resolution mode, though a natural choice of a
translational origin is not possible. Hence a first task is to
`localize' a given probe within $\mathbbm{Z}[\zeta_{n}]$.  To be more
specific, suppose $X$-rays of some planar (quasi)crystalline set $F$
are taken in some directions $o_1,\dots,o_m\in
\mathbbm{Z}[\zeta_{n}]\setminus\{0\}$. Obviously, every point of $F$
is `registered' by every $X$-ray image, hence $F$ is contained in the
\emph{grid}
$$G:=\bigcap_{i=1}^{m}\,\,\left( \bigcup_{v\in F} (v+\mathbbm{R} o_i)\right);$$
see Definition \ref{m-def:thegrid}.  Of course, in general $G$
contains many more points than $F$, hence does not disclose $F$. On
the other hand, only those subsets $F'$ of $G$ whose $X$-rays coincide
with the given data are feasible solutions which lie in a translate of
the underlying model set.  Hence a first problem is to determine the
decomposition of $G$ into the subsets which are compatible with the
underlying $\mathbbm{Z}$-module $\mathbbm{Z}[\zeta_{n}]$, i.e., which
lie in a common translate of $\mathbbm{Z}[\zeta_{n}]$; see Section
\ref{m-sec:sec4} for details.  This problem has its origin in the
practice of quantitative HRTEM since, in general, the $X$-ray
information does not allow us to locate the underlying
$\mathbbm{Z}$-module $\mathbbm{Z}[\zeta_{n}]$.  Using standard results
of algebra, we will actually show much more, namely that the solution
of this \emph{decomposition problem} only depends on $n$ and the given
$\mathbbm{Z}[\zeta_{n}]$-directions but not on the specific $X$-ray
data.  Hence, conceptually, we can consider the different equivalence
classes separately.

Of course, even if a $\mathbbm{Z}[\zeta_{n}]$-equivalence class of the
grid $G$ contains a set $F'$ whose $X$-rays coincide with the given
data, this set need not belong to the underlying cyclotomic model set.
Hence it is clear that additional constraints that are induced by the
construction rules of the underlying model set have to be satisfied to
guarantee feasibility. One possible approach could be to first
reconstruct a potential solution that is compatible with the given
$X$-ray data and then check whether it actually belongs to the
underlying model set. Unfortunately, this approach does not lead to an
efficient algorithm (see Remark \ref{m-rem:firstsep}). Therefore, we
use the specific structure of model sets (originating from some window
through a projection process) and determine which subsets of $G$ can
possibly arise. In fact, all possible solutions (that might actually
lie `far out' in the defining model set) can be found and explored by
translating the given window; see Section \ref{m-sec:sec4}. For many
types of windows, this \emph{separation problem} can be handled by
geometric techniques based on the theory of arrangements; see Section
\ref{m-sec:results}.

\vskip .2cm
The present paper is organized as follows.

As a service to the reader, we begin with two preliminary sections
that put together the notions required and recall several tools from
algebra and the mathematical theory of quasicrystals. In fact, the
algebra is needed not only to properly explain cyclotomic model sets
but is crucial for devising algorithms for checking containment of
points in this structure and also yields best-known bounds for the
running time of our basic algorithms.  (Of course, in view of the
prominent role of group theory in crystallography and materials
science, the relevance of algebraic methods for our cyclotomic
structures does not really come as a surprise.)  In Section
\ref{sec2}, we explain the algebraic concepts in an elementary way
while Section \ref{sec3} gives a concise but sufficiently detailed
account of \emph{model sets}.  In particular, we introduce the special
class of \emph{cyclotomic model sets}, which will be the central
objects of the present paper.  Some examples illustrate the structure
and the beauty of cyclotomic model sets.

The key problems and main results will be formulated in Section
\ref{m-sec:sec4}; their proofs will be given in Section
\ref{m-sec:results}.

\section{Algebraic Background and Notation}
\label{sec2}

For all $n \in \mathbbm{N}$, and $\zeta_{n}$ a fixed primitive $n$th
root of unity in $\mathbbm{C}$ ({\it e.g.}, $\zeta_{n} = e^{2\pi i/n}$), let
$\mathbbm{Q}(\zeta_{n})$ be the corresponding cyclotomic field, i.e.,
the smallest intermediate field of the field extension
$\mathbbm{C}/\mathbbm{Q}$ that contains $\zeta_{n}$. Further, denoting
by $\bar{\zeta}_{n}$ the complex conjugate of $\zeta_{n}$, it is well
known that $\mathbbm{Q}(\zeta_{n}+\bar{\zeta}_{n})$ (defined
analogously) is the maximal real subfield of $\mathbbm{Q}(\zeta_{n})$,
i.e.,
$$\mathbbm{Q}(\zeta_{n})\cap\mathbbm{R}=
\mathbbm{Q}(\zeta_{n}+\bar{\zeta}_{n})\,;$$
see~Washington (1997, p. 15). Throughout this text, we shall use the
notation
$$\mathbbm{K}_{n}=
\mathbbm{Q}(\zeta_{n}),\; \mathbbm{k}_{n}=
\mathbbm{Q}(\zeta_{n}+\bar{\zeta}_{n}),\; \mathcal{O}_{n}=
\mathbbm{Z}[\zeta_{n}],\; 
\thinspace\scriptstyle{\mathcal{O}}
\displaystyle_{n} =\mathbbm{Z}[\zeta_{n}+\bar{\zeta}_{n}]\, ,$$
where $\mathbbm{Z}[\zeta_{n}]$ (resp.,
$\mathbbm{Z}[\zeta_{n}+\bar{\zeta}_{n}]$) is defined as the smallest
subring of $\mathbbm{C}$ that contains $\mathbbm{Z}$ and $\zeta_{n}$
(resp., $\mathbbm{Z}$ and $\zeta_{n}+\bar{\zeta}_{n}$). Further,
$\phi$ will always denote Euler's phi-function (often also called
Euler's totient function), i.e., $$\phi(n) =
\operatorname{card}\left(\{k \in \mathbbm{N}\, |\,1 \leq k \leq n
  \textnormal{ and } \operatorname{gcd}(k,n)=1\}\right)\, .$$
Occasionally, we shall identify $\mathbbm{C}$ with $\mathbbm{R}^{2}$.

The set $\mathcal{O}_{n}$ hosts the corresponding cyclotomic model
sets ({\it cf.} Section \ref{m-sec:cms}); $\mathbbm{K}_{n}$,
$\mathbbm{k}_{n}$, and $\scriptstyle{\mathcal{O}}\displaystyle_{n}$
will be needed for the analysis of the algebraic structure of
$\mathcal{O}_{n}$ that will allow the relevant algorithmic
computations. The following lemma shows how $\mathcal{O}_{n}$ is
related to $\thinspace\scriptstyle{\mathcal{O}}\displaystyle_{n}$.

\begin{lemma}\label{Oo}
For $n \geq 3$, one has: 
\begin{itemize}
\item[{\rm (a)}] $\mathcal{O}_{n}$ is an
  $\thinspace\scriptstyle{\mathcal{O}}\displaystyle_{n}$-module of
  rank $2$.  More precisely, one has $\mathcal{O}_{n} =\,
  \thinspace\scriptstyle{\mathcal{O}}\displaystyle_{n} +\,
  \thinspace\scriptstyle{\mathcal{O}}\displaystyle_{n} \,\zeta_{n}$,
  and $\{1,\zeta_{n}\}$ is an
  $\thinspace\scriptstyle{\mathcal{O}}\displaystyle_{n}$-basis of
  $\mathcal{O}_{n}$.
\item[{\rm (b)}] $\mathbbm{K}_{n}$ is a $\mathbbm{k}_{n}$-vector space
  of dimension $2$.  More precisely, one has $\mathbbm{K}_{n} =\,
  \mathbbm{k}_{n} +\, \mathbbm{k}_{n} \,\zeta_{n}$, and
  $\{1,\zeta_{n}\}$ is a $\mathbbm{k}_{n}$-basis of $\mathbbm{K}_{n}$.
\end{itemize}
\end{lemma}
\begin{proof}
  First, we show (a). The linear independence of $\{1,\zeta_{n}\}$
  over $\thinspace\scriptstyle{\mathcal{O}}\displaystyle_{n}$ is
  clear: by our assumption $n \geq 3$, $\{1,\zeta_{n}\}$ is even
  linearly independent over $\mathbbm{R}$.  For the remainder of the
  assertion we prove that all non-negative integral powers
  $\zeta_{n}^i$ satisfy $\zeta_{n}^i= \alpha + \beta \zeta_{n}$ for
  suitable $\alpha, \beta \in
  \thinspace\scriptstyle{\mathcal{O}}\displaystyle_{n}$.  Using
  induction, it suffices to show $\zeta_{n}^2= \alpha + \beta
  \zeta_{n}$ for suitable $\alpha, \beta \in
  \thinspace\scriptstyle{\mathcal{O}}\displaystyle_{n}$. To this end,
  note that $\bar{\zeta}_{n} = \zeta_{n}^{-1}$ and observe that
  $\zeta_{n}^2 = -1 + (\zeta_{n} + \zeta_{n}^{-1}) \zeta_{n}$.

Claim (b) follows similarly. 
\end{proof}

\begin{figure}
\begin{minipage}{\textwidth}
\includegraphics[scale=.65]{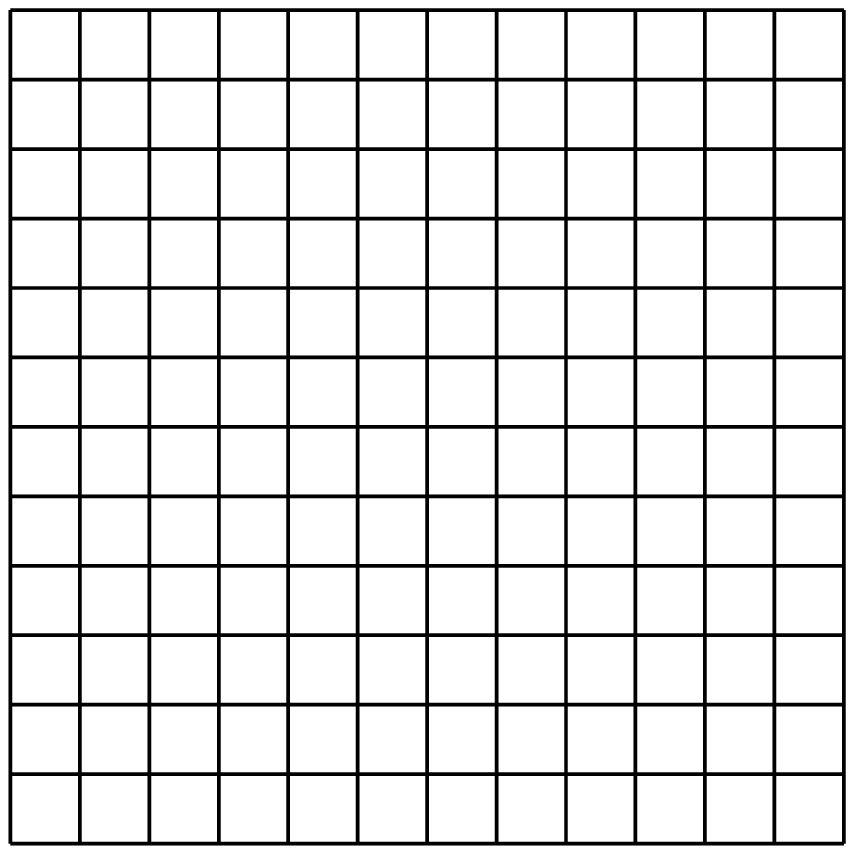}
\hspace{\fill}
\includegraphics[scale=.55]{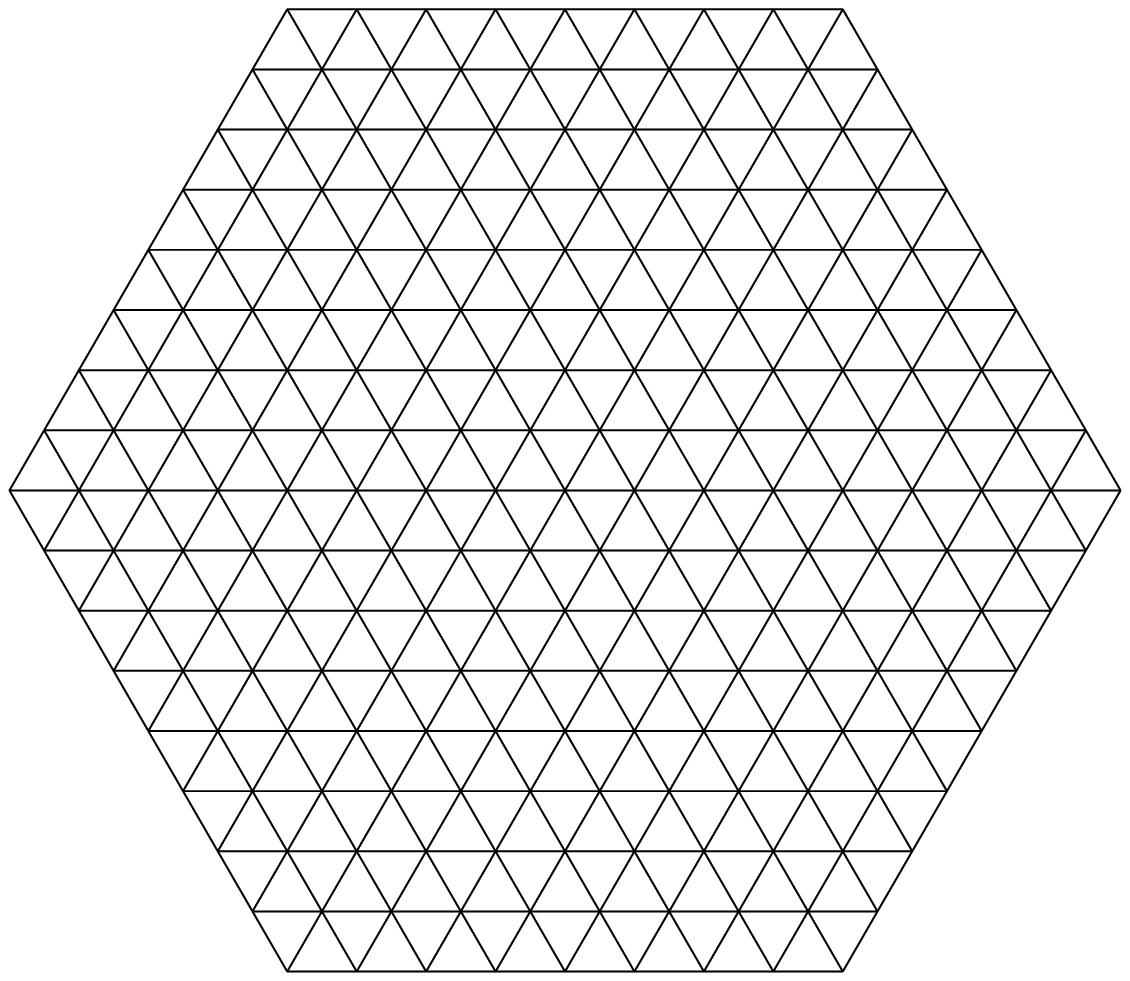}
\label{fig:sqttt}
\end{minipage}
\caption{Central patches of the square tiling (left) and triangular tiling (right).}
\end{figure}

\begin{remark}\label{r0}
  Seen as a point set of $\mathbbm{R}^2$, $\mathcal{O}_{n}$ has
  $N$-fold cyclic symmetry, where
\begin{equation}\label{eq}
N=N(n):= \operatorname{lcm}(n,2)=\left\{
\begin{array}{ll}
n, & \mbox{if $n$ is even,}\\
2n, & \mbox{if $n$ is odd. }
\end{array}\right.
\end{equation}
Except for the one-dimensional case $n\in\{1,2\}$
($\mathcal{O}_{1}=\mathcal{O}_{2}=\mathbbm{Z}$), the crystallographic
cases $n\in\{3,6\}$ (triangular lattice
$\mathcal{O}_{3}=\mathcal{O}_{6}$, see Figure~\ref{fig:sqttt}) and
$n=4$ (square lattice $\mathcal{O}_{4}$, see Figure~\ref{fig:sqttt}),
$\mathcal{O}_{n}$ is dense in $\mathbbm{R}^2$. For the latter, note
that, by Lemma~\ref{Oo}, $\mathcal{O}_{n}$ is an
$\thinspace\scriptstyle{\mathcal{O}}\displaystyle_{n}$-module of rank
$2$, whose $\mathbbm{R}$-span is all of $\mathbbm{R}^2$. For
$n\in\mathbbm{N}\setminus\{1,2,3,4,6\}$,
$\thinspace\scriptstyle{\mathcal{O}}\displaystyle_{n}$ is a
$\mathbbm{Z}$-module of rank $\geq 2$ (see Remark~\ref{r7} below)
embedded in $\mathbbm{R}$, hence a dense set in $\mathbbm{R}$.
Consequently, $\mathcal{O}_{n}$ is then a dense set in
$\mathbbm{R}^2$.
\end{remark}

The following well-known result is needed later to actually compute
the coordinates of $\mathcal{O}_{n}$-points. As usual, $R^{\times}$
denotes the group of units of a given ring $R$.

\begin{proposition}[Gau\ss]\label{gau}
  One has $[\mathbbm{K}_{n} : \mathbbm{Q}] = \phi(n)$ and
  $\{1,\zeta_{n},\zeta_{n}^2, \dots ,\zeta_{n}^{\phi(n)-1}\}$ is a
  $\mathbbm{Q}$-basis of $\mathbbm{K}_{n}$. Moreover, the
  field extension $\mathbbm{K}_{n}/ \mathbbm{Q}$ is a Galois extension
  with Abelian Galois group $G(\mathbbm{K}_{n}/ \mathbbm{Q}) \cong
  (\mathbbm{Z} / n\mathbbm{Z})^{\times}$, where $a\,
  (\textnormal{mod}\, n)$ corresponds to the automorphism given by\/
  $\zeta_{n} \longmapsto \zeta_{n}^{a}$.
\end{proposition}

\begin{proof} See Theorem 2.5 of Washington (1997) and, for the
  statement about the $\mathbbm{Q}$-basis, the proof of Proposition
  1.4 in Chapter V.1 of Lang (1993).
\end{proof}

\begin{remark} 
  Note the identity $$(\mathbbm{Z} / n\mathbbm{Z})^{\times}=\{a\,
  (\textnormal{mod}\, n)\,|\,(a,n)=1\}$$ and consult Table 3 of Baake
  \& Grimm (2004) for examples of the explicit structure of
  $G(\mathbbm{K}_{n}/ \mathbbm{Q})$.
\end{remark} 

\begin{corollary}\label{cr5}
  If $n\geq 3$, one has $[\mathbbm{k}_{n} : \mathbbm{Q}] =
  \phi(n)/2$. Moreover, a $\mathbbm{Q}$-basis of $\mathbbm{k}_{n}$ is
  given by the set $\{1,(\zeta_{n}+\bar{\zeta}_{n}),(\zeta_{n}+\bar{\zeta}_{n})^2,
  \dots ,(\zeta_{n}+\bar{\zeta}_{n})^{\phi(n)/2-1}\}$.
\end{corollary}
\begin{proof}
  The statement about the degree $[\mathbbm{k}_{n} : \mathbbm{Q}]$ is
  an immediate consequence of Lemma~\ref{Oo}(b), Proposition~\ref{gau}
  and the `degree formula' for field extensions: If $E/F/K$ is an
  extension of fields, one has $[E:K] = [E:F] [F:K]$ ({\it cf.} Chapter V.1,
  Proposition 1.2 of Lang (1993)). The statement about the
  $\mathbbm{Q}$-basis again follows from the proof of Proposition 1.4
  in Chapter V.1 of Lang (1993).
\end{proof}

A full $\mathbbm{Z}$-module (i.e., a module of full rank) in an
algebraic number field $\mathbbm{K}$ which contains the number $1$ and
is a ring is called an \emph{order} of $\mathbbm{K}$. It turns out
that among the various orders of $\mathbbm{K}$ there is one
\emph{maximal order} which contains all the other orders, namely the
ring of integers in $\mathbbm{K}$; see Chapter 2, Section 2 of
 Borevich \& Shafarevich (1966). For cyclotomic fields, one has the
following well-known result.

\begin{proposition}\label{p1}
For $n\in \mathbbm{N}$, one has:
\begin{itemize}
\item[{\rm (a)}] $\mathcal{O}_{n}$ is the ring of cyclotomic integers
  in $\mathbbm{K}_{n}$, and hence is its maximal order.
\item[{\rm (b)}]
  $\thinspace\scriptstyle{\mathcal{O}}\displaystyle_{n}$ is the ring
  of integers of $\mathbbm{k}_{n}$, and hence is its maximal order.
\end{itemize}
\end{proposition}
\begin{proof}
See Theorem 2.6 and Proposition 2.16 of Washington (1997).
\end{proof}

\begin{remark}\label{r7} 
  It follows from Proposition~\ref{p1}(a) and Proposition~\ref{gau}
  that $\mathcal{O}_{n}$ is a $\mathbbm{Z}$-module of rank $\phi(n)$
  with $\mathbbm{Z}$-basis $\{1,\zeta_{n},\zeta_{n}^2, \dots
  ,\zeta_{n}^{\phi(n)-1}\}$.  Likewise, Proposition~\ref{p1}(b) and
  Corollary~\ref{cr5} imply that, for $n\geq 3$,
  $\thinspace\scriptstyle{\mathcal{O}}\displaystyle_{n}$ is a
  $\mathbbm{Z}$-module of rank $\phi(n)/2$ with $\mathbbm{Z}$-basis
  given by the set 
  $\{1,(\zeta_{n}+\bar{\zeta}_{n}),(\zeta_{n}+\bar{\zeta}_{n})^2,
  \dots ,(\zeta_{n}+\bar{\zeta}_{n})^{\phi(n)/2-1}\}$.
\end{remark} 

For the subsequent algorithmic computations the minimum polynomial
$\operatorname{Mipo}_{\mathbbm{Q}}(\zeta_{n})$ of $\zeta_{n}$ over
$\mathbbm{Q}$ will be needed since it shows how to replace certain
higher powers of $\zeta_{n}$ by sums of lower ones. As it turns out in
Proposition \ref{gau2}, $\operatorname{Mipo}_{\mathbbm{Q}}(\zeta_{n})$
is simply the following $n$th cyclotomic polynomial.

\begin{definition}\label{fn}
The $n$\emph{th cyclotomic polynomial} is given by 
$$F_{n} := \prod_{\zeta} (X - \zeta)\, ,$$ 
where $\zeta$ runs over all primitive $n$th roots of unity in
$\mathbbm{C}$.
\end{definition}

\begin{lemma}\label{lemma6}
For $n\in \mathbbm{N}$, one has:
\begin{itemize}
\item[{\rm (a)}]
$F_{n}$ is monic and $\operatorname{deg}(F_{n}) = \phi(n)$.
\item[{\rm (b)}]
${{\underset{d \mid n}{\prod} F_{d}}} = X^n - 1$.
\item[{\rm (c)}]
$F_{n} \in \mathbbm{Z}[X]$.
\end{itemize}
\end{lemma}
\begin{proof}
See Chapter VI.3 of Lang (1993).
\end{proof}

\begin{remark}\label{recur}
  Lemma~\ref{lemma6} shows that we can compute the $n$th cyclotomic
  polynomial recursively by use of the Euclidean algorithm in
  $\mathbbm{Z}[X]$.
\end{remark}

\begin{proposition}[Gau\ss]\label{gau2}
The minimum polynomial $\operatorname{Mipo}_{\mathbbm{Q}}(\zeta_{n})$ of 
$\zeta_{n}$ over $\mathbbm{Q}$ is the $n$th cyclotomic polynomial $F_{n}$. 
\end{proposition}

\begin{proof}
  By Definition~\ref{fn}, $\zeta_{n}$ is a root of $F_{n}$. Now, note
  that $\operatorname{Mipo}_{\mathbbm{Q}}(\zeta_{n})$ is, by
  definition, the (uniquely determined) monic polynomial in
  $\mathbbm{Q}[X]$ of minimal degree having $\zeta_{n}$ as a root. Of
  course, it is a standard fact that
  $\operatorname{deg}(\operatorname{Mipo}_{\mathbbm{Q}}(\zeta_{n}))=[\mathbbm{K}_{n}:\mathbbm{Q}]$;
  see Proposition 1.4 in Chapter V.1 of Lang (1993).  By
  Proposition~\ref{gau}, one has
  $[\mathbbm{K}_{n}:\mathbbm{Q}]=\phi(n)$, hence the result follows
  from Lemma~\ref{lemma6}.
\end{proof}

The final result of this preliminary section will provide a uniform
finite upper bound on the number of $\mathcal{O}_{n}$-equivalence
classes in arbitrary grids for given $X$-ray directions.

\begin{proposition}\label{full}
  If $G$ is a torsion-free Abelian group of rank $r$, and $H$ is a
  subgroup which is also of rank $r$, then the subgroup index $[G:H]$
  is finite and equals the absolute value of the determinant of the
  transition matrix $A$ from any\/ $\mathbbm{Z}$-basis of $G$ to any\/
  $\mathbbm{Z}$-basis of $H$.
\end{proposition}

\begin{proof} See Chapter 2, Lemma 6.1.1 of Borevich \& Shafarevich
  (1966).
\end{proof}

\section{Model Sets}\label{sec3}

Now we will first give a brief introduction to model sets and then we
define the class of cyclotomic model sets that will be the underlying
ground structure for the present paper.

\subsection{General Setting}\label{generalsetting}

By definition, \emph{model sets} arise from so-called \emph{cut and
  project schemes}. These are commutative diagrams of the following
form; compare with~Moody (2000) and see Baake {\it et al.} (2002) for a gentle
introduction with many illustrations.
\begin{equation}\label{cutproj0}
  \renewcommand{\arraystretch}{1.4}
\begin{array}{ccccc}
  \,\,\,& \pi & & \pi_{\textnormal{\tiny int}} & \vspace*{-2.0ex} \\
  \mathbbm{R}^{k}\,\,\, & \longleftarrow & \mathbbm{R}^{k}\times H  & \longrightarrow & H \vspace*{1.5ex}\\
  \cup\,\,\,\,&&\,\,\,\,\,\,\,\,\,\,\,\,\,\,\,\cup\mbox{\tiny\, lattice}&&\,\,\,\,\,\,\,\,\,\,\,\,\cup\mbox{\tiny\, dense}\\
  \,\,\,& \mbox{\tiny 1--1} & &  & \vspace*{-2.0ex} \\
  \!\!\!\pi[\widetilde{L}]\,\,\,& \longleftrightarrow &\, \widetilde{L}& \longrightarrow &\pi_{\textnormal{\tiny
      int}}[\widetilde{L}] \\
\end{array}
\end{equation}
Here, $H$ is some locally compact Abelian group, $\pi$ and
$\pi_{\textnormal{\tiny int}}$ are the canonical projections, and
$\widetilde{L}$ is a lattice in $\mathbbm{R}^{k}\times H$, i.e.,
$\widetilde{L}$ is a discrete subgroup of $\mathbbm{R}^{k}\times H$
such that the quotient group
$$(\mathbbm{R}^{k}\times H)\,\big/\,\widetilde{L}$$ is compact. Further,
$\pi_{\textnormal{\tiny int}}[\widetilde{L}]$ is a dense subset of $H$
and the restriction of $\pi$ to $\widetilde{L}$ is assumed to be
injective. Writing $L:=\pi[\widetilde{L}]$, one can define a map
$.^{\star}\! : \, L \longrightarrow H$ by $x\longmapsto
\pi_{\textnormal{\tiny int}}(\pi|_{L}^{-1}(x))$. Then, one has
$[L]^{\star}=\pi_{\textnormal{\tiny int}}[\widetilde{L}]$. If the map
$.^{\star}$ is injective, we denote the inverse of its co-restriction
$.^{\star}\! : \, L \longrightarrow [L]^{\star}$ by $.^{-\star } \! :
\, [L]^{\star} \longrightarrow L$.

In the following we use the notation $A^{\circ}$, $\overline{A}$,
$\partial A$ for the standard topological operators \emph{interior},
\emph{closure}, and \emph{boundary} of a set $A$ in a locally compact
Abelian group.
  
\begin{definition}\label{m-def:modelsets}
\begin{itemize}
\item[{\rm (a)}] Given the cut and project scheme $(\ref{cutproj0})$,
  a subset $W\subset H$ is called a \emph{window} if $\varnothing\,
  \neq\, W^{\circ}\subset W\subset \overline{W^{\circ}}$ and
  $\overline{W^{\circ}}$ is compact.
\item[{\rm (b)}] Given any window $W\subset H$, and any
  $t\in\mathbbm{R}^d$, we obtain a \emph{model set}
$$\varLambda(t,W) := t+\varLambda(W)$$
relative to the cut and project scheme by setting
$$\varLambda(W):=\{x\in L\,|\,x^{\star}\in W\}\, .$$
Further, $\mathbbm{R}^{k}$ $($resp., $H$$)$ is called the
\emph{physical} $($resp., \emph{internal}$)$ space and $W$ is also
referred to as the \emph{window} of $\varLambda(t,W)$. The map
$.^{\star}\! : \, L \longrightarrow H$, as defined above, is the
so-called \emph{star map}.
\end{itemize}
\end{definition}

For details about model sets and general background material see
 Moody (2000) and Baake \& Moody (2000); see~Baake {\it et al.} (2002) for
detailed graphical illustrations of the projection method.

\begin{remark}\label{transl}
  The translation vector $t$ in Definition \ref{m-def:modelsets}
  stresses an intrinsic character of model sets. While the structure
  model specifies the cut and project scheme $k$, $H$ and
  $\widetilde{L}$, and also the window $W$, a natural choice of the
  origin is usually not possible.
\end{remark}

\begin{remark}
  Without loss of generality, we may assume that the stabilizer
  $H_{W}$ of the window $W$, i.e.,
$$
H_{W}:=\{h\in H\,|\,h+W=W\},
$$
is the trivial subgroup of $H$, i.e., $H_{W}=\{0\}$.  Observe that the
latter is always the case if $H$ is some Euclidean space, i.e., if one
has $H=\mathbbm{R}^d$ for some suitable $d\in\mathbbm{N}$. Note
further that the star map is a homomorphism of Abelian groups.
\end{remark}

The following remark collects some properties of model sets; for
details see Moody (2000).

\begin{remark}\label{propms1}
  In the following, for $x\in\mathbbm{R}^d$ and $r>0$, we denote by
  $B_{r}(x)$ the open ball of radius $r$ about $x$. The model set
  $\varLambda:=\varLambda(t,W)\subset \mathbbm{R}^d$ is a \emph{Delone
    set}, meaning that $\varLambda$ is both uniformly discrete (i.e.,
  there is a radius $r>0$ such that every ball of the form $B_{r}(x)$,
  where $x\in\mathbbm{R}^{d}$, contains at most one point of
  $\varLambda$) and relatively dense (i.e., there is a radius $R>0$
  such that every ball of the form $B_{R}(x)$, where
  $x\in\mathbbm{R}^{d}$, contains at least one point of $\varLambda$).
  Also, $\varLambda$ has \emph{finite local complexity}, i.e.,
  $\varLambda-\varLambda$ is discrete and closed.  (Note that
  $\varLambda$ has finite local complexity iff for every $r>0$ there
  are, up to translation, only finitely many point sets (called
  \emph{patches of diameter} $r$) of the form $\varLambda\cap
  B_{r}(x)$, where $x\in\mathbbm{R}^d$.) In fact, $\varLambda$ is even
  a \emph{Meyer set} (i.e., in addition, $\varLambda-\varLambda$ is
  uniformly discrete).

  Further, $\varLambda$ is \emph{aperiodic}, i.e., has no
  translational symmetries iff the star map is injective. In fact, the
  kernel of the star map is the group of translational symmetries of
  $\varLambda$.

  If $\varLambda$ is \emph{regular}, i.e., the boundary $\partial W$
  of the window $W$ has (Haar) measure $0$ in $H$, then $\varLambda$
  is \emph{pure point diffractive} ({\it cf.}~Schlottmann, 2000).  If
  $\varLambda$ is \emph{generic}, i.e., $[L]^{\star}\cap\, \partial W
  = \varnothing$, then $\varLambda$ is \emph{repetitive}. This means
  that, given any patch of radius $r$, there is a radius $R$ such that
  any ball $B_{R}(x)$ in $\mathbbm{R}^d$ contains at least one
  translate of this patch; see~Schlottmann (2000). If $\varLambda$ is
  both generic and regular, the frequency of repetition of finite
  patches is well defined, i.e., for every finite patch, the number of
  occurrences of translates of this patch per unit volume in the ball
  $B_{r}(0)$ of radius $r$ about the origin $0$ approaches a positive
  limit as $r\rightarrow \infty$; {\it cf.}~Schlottmann (1998).
\end{remark}

For the discrete tomography of aperiodic model sets, one additional
difficulty, in comparison to the crystallographic case, stems from the
fact that it is not sufficient to consider one pattern and its
translates to define the setting. In particular, to define the
analogue of a specific crystal, one has to add all infinite patterns
that emerge as limits of sequences of translates defined in the local
topology (LT). Here, two patterns are $\varepsilon$-close if, after a
translation by a distance of at most $\varepsilon$, they agree on a
ball of radius $1/\varepsilon$ around the origin. If the starting
pattern $P$ is crystallographic, no new patterns are added; but if $P$
is a generic aperiodic model set, one ends up with uncountably many
different patterns, even up to translations! Nevertheless, all of them
are locally indistinguishable (LI).  This means that every
\emph{finite} patch in $\varLambda$ also appears in any of the other
elements of the LI-class and vice versa; see~Baake (2002) for
details.

\begin{remark}\label{m-rem:LI}
  The entire LI-class of a regular, generic model set $\varLambda(W)$
  can be shown to consist of all sets $t+\varLambda(\tau+W)$, with $t
  \in \mathbbm{R}^{d}$ and $\tau$ such that $[L]^{\star}\cap
  \partial(\tau+W) =\varnothing$ (i.e., $\tau$ is in a generic
  position), and all patterns obtained as limits of sequences
  $t+\varLambda(\tau_{n}+W)$, with all $\tau_{n}$ in a generic
  position; see~Baake (2002). Each such limit is then a \emph{subset}
  of some $t+\varLambda(\tau+W)$, as $\tau$ might not be in a generic
  position. In view of this complication, we must make sure that we
  reconstruct finite subsets of \emph{generic} model sets. This will
  be reflected in Definitions~\ref{wsets}~and~\ref {m-def:consrecprob}
  of Section \ref{m-sec:sec4}.
\end{remark}

\subsection{Cyclotomic Model Sets}\label{m-sec:cms}

In the present paper we will study the discrete tomography of a
special class of planar model sets, the \emph{cyclotomic} model sets,
which can be described in algebraic terms and have an Euclidean
internal space. In the following let $n\in \mathbbm{N} \setminus
\{1,2\}$.

Before we formally introduce the cut and project scheme from which the
\emph{cyclotomic} model sets arise, let us consider some main
ingredients.

The elements of the Galois group $G(\mathbbm{K}_{n}/ \mathbbm{Q})$
(see Proposition~\ref{gau}) come in pairs of complex conjugate
automorphisms. Let the set $\{\sigma_{1},\dots,\sigma_{\phi(n)/2}\}$
arise from $G(\mathbbm{K}_{n}/ \mathbbm{Q})$ by choosing exactly one
automorphism from each such pair. Here, we always choose $\sigma_{1}$
as the identity rather than the complex conjugation. Every such choice
induces a map
$$
\,\,.\,\widetilde{\hphantom{a}}\,:\, \mathcal{O}_{n}\longrightarrow
(\mathbbm{R}^2)^{\frac{\phi(n)}{2}}
$$
through
$$
z\longmapsto \left(z,\sigma_{2}(z),\dots,\sigma_{\frac{\phi(n)}{2}}(z)\right)\,.
$$
(Actually, $.\,\widetilde{\hphantom{a}}$ and the following map
$.^{\star}$ are defined on $\mathbbm{K}_{n}$, but it is their
restriction to $\mathcal{O}_{n}$ that is relevant here.)

With the understanding that for $\phi(n)=2$ (i.e., $n\in\{3,4,6\}$),
the singleton
$$(\mathbbm{R}^2)^{\frac{\phi(n)}{2}-1}=(\mathbbm{R}^2)^{0}$$
is the trivial (locally compact) Abelian group $\{0\}$ each such
choice induces a map
$$.^{\star}\! : \, \mathcal{O}_{n} \longrightarrow
(\mathbbm{R}^2)^{\frac{\phi(n)}{2}-1}\,,$$ defined by
$.^{\star}\!:\equiv 0$, if $n\in\{3,4,6\}$, and
$$z\longmapsto
\left(\sigma_2(z),\dots,\sigma_{\frac{\phi(n)}{2}}(z)\right)$$
otherwise.  Then, $[\mathcal{O}_{n}]\widetilde{\hphantom{a}}$ is a
Minkowski representation of the maximal order $\mathcal{O}_{n}$ of
$\mathbbm{K}_{n}$, see~Chapter 2, Section 3 of Borevich \&
Shafarevich (1966) and Theorem 2.6 of Washington (1997). It follows
that $[\mathcal{O}_{n}]\widetilde{\hphantom{a}}$ is a (full) lattice
in $\mathbbm{R}^{2}\times(\mathbbm{R}^2)^{\phi(n)/2-1}$. Here, since
the space 
$\mathbbm{R}^{2}\times(\mathbbm{R}^2)^{\phi(n)/2-1}$ is Euclidean,
this means that there are $\phi(n)$ $\mathbbm{R}$-linearly independent
vectors in $\mathbbm{R}^{2}\times(\mathbbm{R}^2)^{\phi(n)/2-1}$ having
the property that $[\mathcal{O}_{n}]\widetilde{\hphantom{a}}$ is the
$\mathbbm{Z}$-span of these vectors; compare Chapter 2, Sections 3 and
4 of Borevich \& Shafarevich (1966). In fact, the set
$$\left\{1\widetilde{\hphantom{a}},(\zeta_{n})\widetilde{\hphantom{a}},\dots,
  (\zeta_{n}^{\phi(n)-1})\widetilde{\hphantom{a}}\right\}$$ has this
property; {\it cf.} Proposition~\ref{p1} and Remark~\ref{r7}.  Further, the
image $[\mathcal{O}_{n}]^{\star}$ is dense in
$(\mathbbm{R}^2)^{\phi(n)/2-1}$. This follows for instance from the
existence of a Pisot number of (full) degree $\phi(n)/2$ in
$\thinspace\scriptstyle{\mathcal{O}}\displaystyle_{n}$; see
Chapter~$2$, Section~$3$ of Borevich \& Shafarevich (1966) and
 Pleasants (2000). Multiplication by such a Pisot number in the
physical space then translates via the map $.^{\star}$ into a
contraction in all directions of the internal space, as defined by the
$\mathbbm{Q}$-span of the projected basis vectors of the lattice.

Now, the cyclotomic model sets arise from cut and project schemes of
the following form, where we follow Moody (2000), modified in
the spirit of the algebraic setting of Pleasants (2000).

\begin{equation}\label{cutproj1}
\renewcommand{\arraystretch}{1.4}
\begin{array}{ccccc}
& \pi & & \pi_{\textnormal{\tiny int}}^{} & \vspace*{-2.0ex} \\
\!\!\!\!\!\mathbbm{R}^{2} & \longleftarrow & \;\;\;\;\mathbbm{R}^{2}\times(\mathbbm{R}^2)^{\frac{\phi(n)}{2}-1}  & 
\longrightarrow & (\mathbbm{R}^2)^{\frac{\phi(n)}{2}-1}\vspace*{1.8ex} \\
\!\!\!\!\!\cup\mbox{\tiny\,  }&&\,\,\cup\mbox{\tiny\, lattice}&&\!\cup\mbox{\tiny\, dense}\\
 & \mbox{\tiny 1--1} & &  & \vspace*{-2.0ex} \\
\!\!\!\!\! \mathcal{O}_{n}& \longleftrightarrow &\!\!\!\!\!\!\!\! [\mathcal{O}_{n}]\widetilde{\hphantom{a}}
& 
\longrightarrow &\!\!\!\!\!\!\!\!\!\!\![\mathcal{O}_{n}]^{\star} \\
\end{array}
\end{equation}

As described above, one has
$$[\mathcal{O}_{n}]\widetilde{\hphantom{a}}=\Big\{\big(z,\underbrace{(\sigma_{2}(z),\dots ,\sigma_{\frac{\phi(n)}{2}}(z))}_{=z^{\star}}\big)\,\Big | \, 
z\in \mathcal{O}_{n}\Big\}\,.$$

Recall that for $n\ne 3,4,6$ also the first inclusion
$\mathcal{O}_{n}\subset \mathbbm{R}^{2}$ involves a dense set. Now
here is the definition of the class of cyclotomic model sets; for more
details and related general algebraic settings, see Pleasants (2000).
 
\begin{definition}\label{defcyc}
  Given any window $W\subset(\mathbbm{R}^2)^{\phi(n)/2-1}$, and any
  $t\in\mathbbm{R}^2$, we obtain a planar model set
  $$\varLambda_{n}(t,W) := t+\varLambda_{n}(W)$$ relative to the above
  cut and project scheme~$(\ref{cutproj1})$ $($i.e., relative to any
  choice of the set\linebreak $\{\sigma_{i} \,|\, i\in\{2,\dots,\phi(n)/2\}\}$
  as described above$)$ by setting
  $$\varLambda_{n}(W):=\{z\in\mathcal{O}_{n}\,|\,z^{\star}\in W\}\,
  .$$ We set
$$\mathcal{M}(\mathcal{O}_{n}):=
\left\{\varLambda_{n}(t,W)   ~\left | ~
    \begin{array}{l} t\in\mathbbm{R}^2\, ,\,
      W\subset(\mathbbm{R}^2)^{\frac{\phi(n)}{2}-1} \mbox{ is}\\
      \mbox{a window}\end{array} \right \}. \right. $$ Then, the class
$\mathcal{CM}$ of \emph{cyclotomic} model sets is defined as
$$
\mathcal{CM}:=\bigcup_{n\,\in\,\mathbbm{N}\setminus \{1,2\}}\mathcal{M}(\mathcal{O}_{n})\,.
$$
\end{definition}

\begin{remark}
  The set $\varLambda:=\varLambda_{n}(t,W)\subset \mathbbm{R}^2$ is
  aperiodic iff $n\notin\{3,4,6\}$, i.e., the translates of the square
  (resp., triangular) lattice are the only cyclotomic model sets
  having translational symmetries; compare Remark~\ref{propms1}.  If,
  for a given $n$, $\varLambda$ is both generic and regular, and, if
  the window $W$ has $m$-fold cyclic symmetry with $m$ a divisor of
  $\operatorname{lcm}(n,2)$ and all in a suitable representation of
  the cyclic group $\mathsf{C}_m$ of order $m$, then $\varLambda$ has
  $m$-fold cyclic symmetry in the sense of symmetries of LI-classes.
  This means that $\varLambda$ and the structure obtained by applying
  an appropriate `symmetry' are locally indistinguishable (LI); see
   Baake (2002) for details on the symmetry concept.
\end{remark}

\subsubsection{Some Examples}\label{m-sec:examples}

We give five examples of cyclotomic model sets. The first two are
periodic of the form $\varLambda_{n}(0,W)\in
\mathcal{M}(\mathcal{O}_{n})$ with $n\,\in \,\{3,4\}$ (and hence
$W=\{0\}$), while the last three are aperiodic cyclotomic model sets
of the form $\varLambda_{n}(0,W)\in \mathcal{M}(\mathcal{O}_{n})$,
with $n\in\{5,8,12\}$ (whence having an internal space of dimension
$2$).

\begin{itemize}
\item[(a)] The planar, generic, regular and periodic cyclotomic model
  set with $4$-fold cyclic symmetry associated with the well-known
  square tiling is the square lattice, which can be described in
  algebraic terms as $\varLambda_{\text{SQ}}
  :=\varLambda_{4}(0,W)=\mathbbm{Z}[i] = \mathcal{O}_{4}$; see
  Figure~\ref{fig:sqttt}.

\item[(b)] The planar, generic, regular and periodic cyclotomic model
  set with $6$-fold cyclic symmetry associated with the well-known
  triangle tiling is the triangle lattice, which can be described in
  algebraic terms as $\varLambda_{\text{TRI}}
  :=\varLambda_{3}(0,W)=\mathcal{O}_{3}$; see Figure~\ref{fig:sqttt}.

\begin{figure}
\centerline{\epsfysize=0.40\textheight\epsfbox{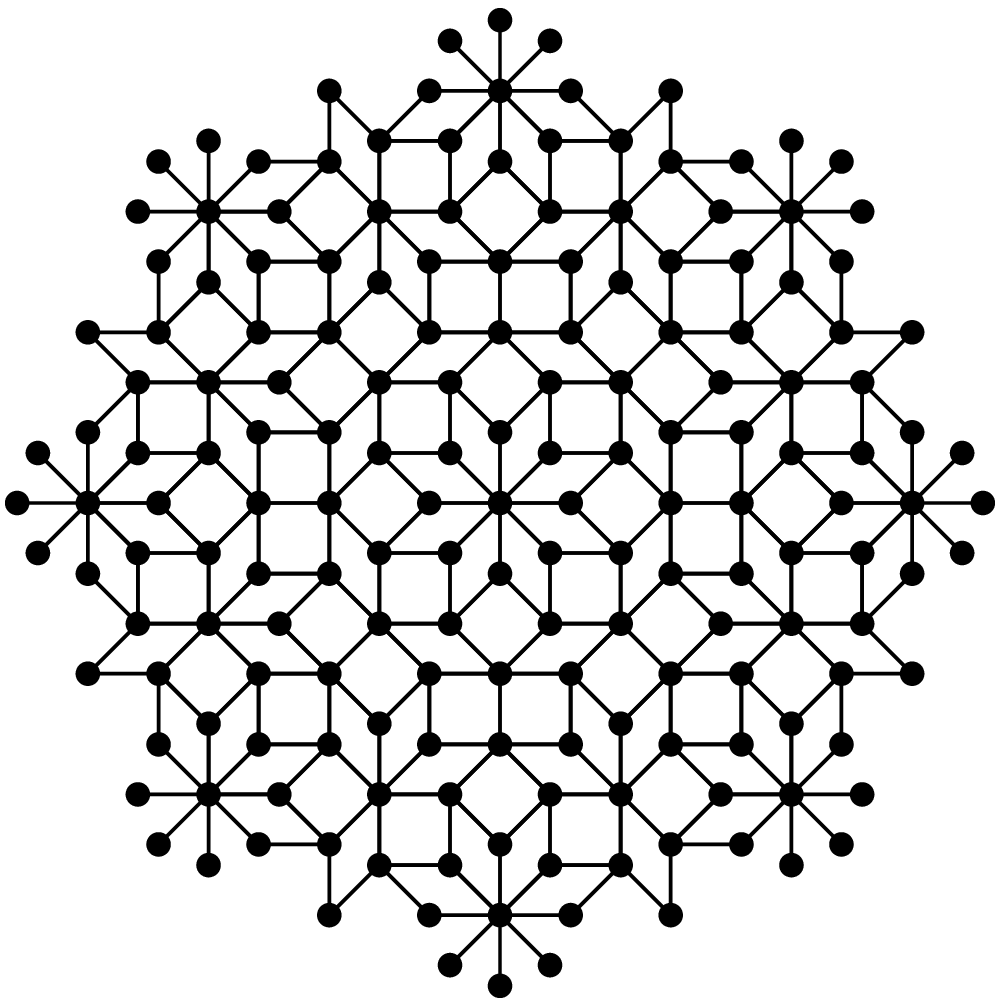}\hspace{0.15\textwidth}
\epsfysize=0.40\textheight\epsfbox{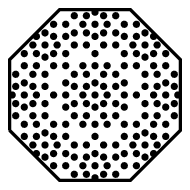}}
\caption{A central patch of the eightfold symmetric Ammann-Beenker tiling with vertex set $\varLambda_{\rm AB}$ (left) and the $.^{\star}$-image of $\varLambda_{\rm AB}$ inside the octagonal window in the internal space (right), with relative scale as described in the text.}
\label{fig:ab}
\end{figure}

\item[(c)] The planar, generic and regular model set with $8$-fold
  cyclic symmetry associated with the Ammann-Beenker tiling~(Baake \&
  Joseph, 1990; Ammann {\it et
  al.}, 1992; G\"ahler, 1993) can be described
  in algebraic terms as
$$\varLambda_{\rm AB} := \{z \in \mathcal{O}_{8}\, | \,z^{\star} \in W\}\, ,$$
where the star map $.^{\star}$ is the Galois automorphism in
$G(\mathbbm{K}_{8}/ \mathbbm{Q})$, defined by $\zeta_{8} \longmapsto
\zeta_{8}^3$, and the window $W$ is the regular octagon centred at the
origin and of unit edge length, with orientation as in
Figure~\ref{fig:ab}.  This construction also gives a tiling with
squares and rhombi, both having edge length $1$; see
Figure~\ref{fig:ab}.

If $t\in\mathbbm{R}^2\setminus\{0\}$ is chosen such that $t+W$ is
again in a generic position (this is true for almost all
$t\in\mathbbm{R}^2$), the replacement of $W$ by $t+W$ again leads to
an Ammann-Beenker tiling. Moreover, the two tilings are locally
indistinguishable (compare Remark \ref{m-rem:LI}).

\item[(d)] The planar and regular model set with $10$-fold cyclic
  symmetry associated with the T\"ubingen triangle tiling~(Baake {\it et
  al.}, 1990a, b) can be described in algebraic terms
  as
$$\varLambda_{\rm TTT}^{t} := \{z \in \mathcal{O}_{5}\, | \,z^{\star} \in t+W\}\, ,$$
where the star map $.^{\star}$ is the Galois automorphism in
$G(\mathbbm{K}_{5}/ \mathbbm{Q})$, defined by $\zeta_{5} \longmapsto
\zeta_{5}^2$. Moreover, the window $W$ is the regular decagon centred
at the origin, with vertices in the directions that arise from the
$10$th roots of unity by a rotation through $\pi/10$, and of edge
length $\tau/\sqrt{\tau+2}$, where $\tau$ is the golden ratio, i.e.,
$\tau=(\sqrt{5}+1)/2$.  Furthermore, $t$ is an element of
$\mathbbm{R}^2$. Note that $\varLambda_{\text{TTT}}^{0}$ is not
generic, while generic examples are obtained for almost all $t\in
\mathbbm{R}^{2}$.  Generic $\varLambda_{\text{TTT}}^{t}$ always give a
triangle tiling with long (short) edges of lengths $1$ ($1/\tau$,
respectively); see Figure~\ref{fig:ttt}.  Different generic choices of
$t$ result in locally indistinguishable (LI) T\"ubingen triangle
tilings (compare again Remark \ref{m-rem:LI}).

\begin{figure}
\centerline{\epsfxsize=0.55\textwidth\epsfbox{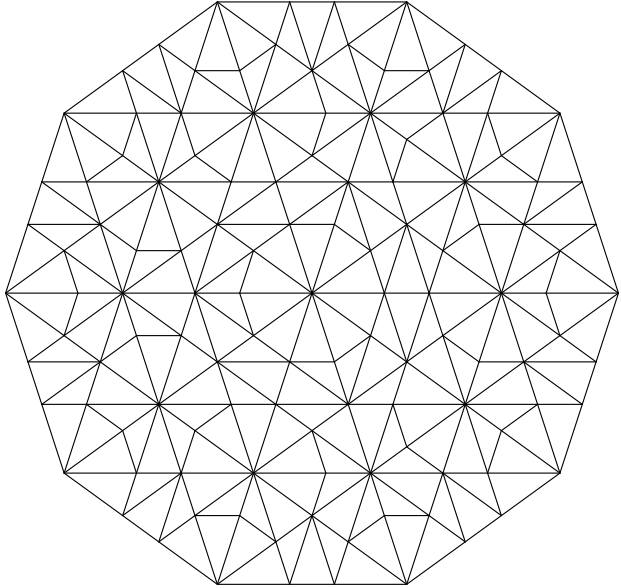}}
\caption{A central patch of the tenfold symmetric T\"ubingen triangle tiling.}
\label{fig:ttt}
\end{figure}

\item[(e)] The planar and regular model set with $12$-fold cyclic
  symmetry associated with the shield tiling~(G\"ahler, 1993) can be
  described in algebraic terms as
$$\varLambda_{\rm S}^{t} := \{z \in \mathcal{O}_{12}\, | \,z^{\star} \in t+W\}\, ,$$
where the star map $.^{\star}$ is the Galois automorphism in
$G(\mathbbm{K}_{12}/ \mathbbm{Q})$, defined by $\zeta_{12} \longmapsto
\zeta_{12}^5$, and the window $W$ is the regular dodecagon centred at
the origin, with vertices in the directions that arise from the $12$th
roots of unity by a rotation through $\pi/12$, and of edge length $1$.
Again, $t$ is an element of $\mathbbm{R}^2$. Note that
$\varLambda_{\text{S}}^{0}$ is not generic, while
$\varLambda_{\text{S}}^{t}$ is generic for almost all $t\in
\mathbbm{R}^{2}$. The shortest distance between points in a generic
$\varLambda_{\text{S}}^{t}$ is $(\sqrt{3}-1)/\sqrt{2}$.  Joining such
points by edges results in a shield tiling, i.e., a tiling with
triangles, squares and so-called shields, all having edge length
$(\sqrt{3}-1)/\sqrt{2}$; see Figure~\ref{fig:s} for a generic example.
Different generic choices of $t$ result in locally indistinguishable
shield tilings (compare again Remark \ref{m-rem:LI}).

\begin{figure}
\centerline{\epsfxsize=0.55\textwidth\epsfbox{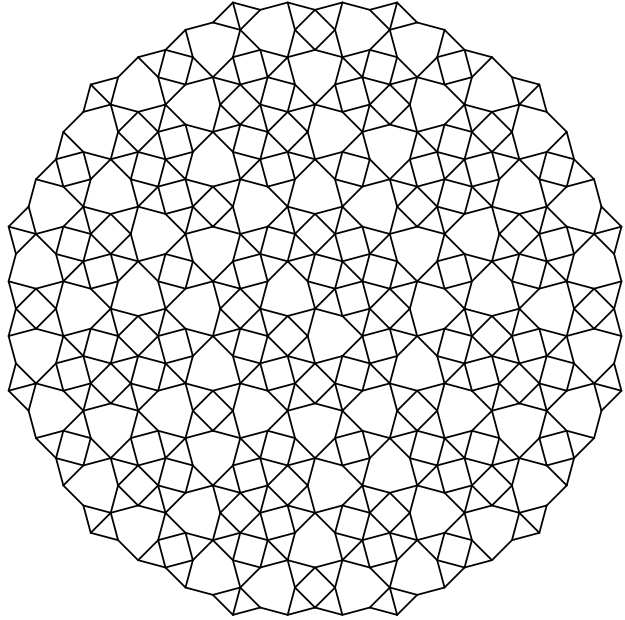}}
\caption{A central patch of the twelvefold symmetric shield tiling.}
\label{fig:s}
\end{figure}

\end{itemize}

\section{Discrete Tomography of Planar Model Sets: Problems and Main Results} \label{m-sec:sec4}

\subsection{Consistency, Reconstruction and Uniqueness}

It is clear that each subset of the lattice $\mathbbm{Z}^2$ is
determined uniquely by one $X$-ray in an irrational direction.
Therefore, the nontrivial classical problems of discrete tomography
involve lattice directions, i.e. directions spanned by two lattice
points.  One now needs the correct analogue of lattice directions in
the framework of cyclotomic model sets.

\begin{definition}
Let $n\in\mathbbm{N}\setminus\{1,2\}$.
\begin{itemize}
\item[\rm{(a)}] The elements of $\mathcal{O}_{n} \setminus\{0\}$ are
  called {\em $\mathcal{O}_{n}$-directions}.
\item[\rm{(b)}] For an $\mathcal{O}_{n}$-direction $o$, we denote by
  $\mathcal{L}_{o}$ the set of lines $t+\mathbbm{R} o$ with $t\in
  \mathbbm{R}$, while $\mathcal{L}^{\mathcal{O}_{n}}_{o}\subset
  \mathcal{L}_{o}$ is the set of module lines in direction $o$, i.e.,
  the set of lines $t+\mathbbm{R} o$ in $\mathbbm{R}^2$ with $t\in
  \mathcal{O}_{n}$.
\end{itemize}
\end{definition}
 
Since every $\mathcal{O}_{n}$-direction is parallel to a non-zero
element of the difference set\linebreak $\varLambda_{n}(t,W) -
\varLambda_{n}(t,W)\subset \mathcal{O}_{n}$ (Huck, 2006), the notion
of $\mathcal{O}_{n}$-directions is indeed the natural extension for
cyclotomic model sets.

\begin{definition}
  Let $n\in\mathbbm{N}\setminus\{1,2\}$ and let $F\subset
  \mathbbm{R}^2$ be a finite set which lives on $\mathcal{O}_{n}$,
  i.e., $F\subset t+\mathcal{O}_{n}$, where $t\in\mathbbm{R}^{2}$.
  Furthermore, let $o$ be an $\mathcal{O}_{n}$-direction.  Then, the
  $($\emph{discrete parallel}\/$)$ \emph{$X$-ray} of $F$ in direction
  $o$ is the function
$$X_{o}F: \mathcal{L}_{o} \longrightarrow \mathbbm{N}_{0}:=\mathbbm{N}\cup\{0\}\,,$$ 
defined by $$X_{o}F(\ell) := \operatorname{card}(F \cap \ell\,)\,.$$
\end{definition}

\begin{remark}
  Obviously, $X_{o}F$ has finite support $\operatorname{supp}(X_{o}F)$
  (the set of lines in direction $o$ that pass through at least one
  point of $F$) and, moreover,
$$\sum_{\ell \in \operatorname{supp}(X_{o}F)}X_{o}F(\ell) = \operatorname{card}(F)\,.$$ 
\end{remark}

In view of the complications with limits indicated at the end of
Section \ref{generalsetting}, we will make sure that we reconstruct
finite subsets of \emph{generic} model sets, i.e., subsets whose
$.^{\star}$-image lies in the \emph{interior} of the window. This
restriction to the generic case is the proper analogue of the
restriction to {\em perfect} lattices and their translates in the
classical case.

\begin{definition}\label{wsets}
  Let $n\in\mathbbm{N}\setminus\{1,2\}$, let $W\subset
  (\mathbbm{R}^2)^{\phi(n)/2-1}$ be a window $(${\it cf.}
  Definition~$\ref{m-def:modelsets}$$)$, and let a star map
  $.^{\star}$ be given, i.e., a map $.^{\star}\! : \,
  \mathcal{O}_{n}\longrightarrow (\mathbbm{R}^2)^{\phi(n)/2-1}$, given
  by $z\longmapsto 0$, if $n\in\{3,4,6\}$, and given by $z\longmapsto
  (\sigma_2(z),\dots,\sigma_{\phi(n)/2}(z))$ otherwise $($as described
  in Definition~$\ref{defcyc}$$)$. Then, the elements of the subset
$$
\{\varLambda_{n}(t,\tau + W^{\circ})\,|\,t\in\mathbbm{R}^{2},\tau \in (\mathbbm{R}^2)^{\frac{\phi(n)}{2}-1} \} 
$$ 
of $\mathcal{M}(\mathcal{O}_{n})$, which are defined by use of the
above star map $.^{\star}$, are called
$W^{\circ}_{\mathcal{M}(\mathcal{O}_{n}),\star}$\emph{-sets}.
\end{definition}

\begin{remark}
  Let $n\in\mathbbm{N}\setminus\{1,2\}$. Note that, if $W\subset
  (\mathbbm{R}^2)^{\phi(n)/2-1}$ is a window, then its interior
  $W^{\circ}$ is also a window. Note further that for $n=4$ (resp.,
  $n\in\{3,6\}$) the set of
  $W^{\circ}_{\mathcal{M}(\mathcal{O}_{n}),\star}$-sets simply
  consists of all translates of the square lattice $\mathcal{O}_{4}$
  (resp., triangular lattice $\mathcal{O}_{3}$).
\end{remark}

\begin{definition}[Consistency, Reconstruction, and Uniqueness
  Problem]\label{m-def:consrecprob}
  Let the data be given as in Definition \ref{wsets}. Further, let
  $o_1,\dots,o_m$ be $m\geq 2$ pairwise non-parallel
  $\mathcal{O}_{n}$-directions. The corresponding consistency,
  reconstruction and uniqueness problems are defined as follows.

  \begin{quote}
    {\sc Consistency}. \\
    Given functions $p_{o_{i}} : \mathcal{L}_{o_{i}} \longrightarrow
    \mathbbm{N}_{0}$, $i\in\{1,\dots,m\}$, whose supports are finite
    and satisfy $\operatorname{supp}(p_{o_{i}})\subset
    \mathcal{L}^{\mathcal{O}_{n}}_{o_{i}}$, decide whether there is a
    finite set $F$ which is contained in a
    $W^{\circ}_{\mathcal{M}(\mathcal{O}_{n}),\star}$-set and satisfies
    $X_{o_{i}}F=p_{o_{i}}$, $i\in\{1,\dots,m\}$.
\end{quote}

\begin{quote}
  {\sc Reconstruction}. \\
  Given functions $p_{o_{i}} : \mathcal{L}_{o_{i}} \longrightarrow
  \mathbbm{N}_{0}$, $i\in\{1,\dots,m\}$, whose supports are finite and
  satisfy $\operatorname{supp}(p_{o_{i}})\subset
  \mathcal{L}^{\mathcal{O}_{n}}_{o_{i}}$, decide whether there exists
  a finite set $F$ in a
  $W^{\circ}_{\mathcal{M}(\mathcal{O}_{n}),\star}$-set that satisfies
  $X_{o_{i}}F=p_{o_{i}}$, $i\in\{1,\dots,m\}$, and, if so, construct
  one such $F$.
\end{quote}

\begin{quote}
  {\sc Uniqueness}. \\
  Given a finite subset $F$ of a
  $W^{\circ}_{\mathcal{M}(\mathcal{O}_{n}),\star}$-set, decide whether
  there is a different finite set $F'$ that is also a subset of a
  $W^{\circ}_{\mathcal{M}(\mathcal{O}_{n}),\star}$-set and satisfies
  $X_{o_{i}}F=X_{o_{i}}F'$, $i\in\{1,\dots,m\}$.
\end{quote}

\end{definition}
Note that the parameter $n$, the directions $o_i$, and the window $W$
are assumed to be fixed, i.e., are \emph{not} part of the input.

For results on the computational complexity of these problems in the
lattice case (and the Turing machine as the model of computation), see
 Gritzmann (1997) and Gardner {\it et al.} (1999).

\subsection{The Decomposition Problem}\label{decomp}

Now we introduce the problem of how to decompose a grid ({\it cf.}
Definition \ref{m-def:thegrid}) into translates of maximal
$\mathcal{O}_{n}$-subsets.  Note that the crystallographic cases,
namely, the triangular lattice and the square lattice, are included.

\begin{definition}\label{m-def:thegrid}
  Let $n\in\mathbbm{N}\setminus\{1,2\}$ and let $o_1,\dots,o_m$ be
  $m\geq 2$ pairwise non-parallel $\mathcal{O}_{n}$-directions.
  Moreover, let $p_{o_{i}} : \mathcal{L}_{o_{i}} \longrightarrow
  \mathbbm{N}_{0}$, $i\in\{1,\dots,m\}$, be functions whose supports
  are finite and satisfy
$$\operatorname{supp}(p_{o_{i}})\subset
\mathcal{L}^{\mathcal{O}_{n}}_{o_{i}}\,.$$ Then, the associated
\emph{grid} $G_{\{p_{o_{i}}|i\in\{1,\dots,m\}\}}$ is defined by
$$G_{\{p_{o_{i}}|i\in\{1,\dots,m\}\}}:=\bigcap_{i=1}^{m}\,\,\left( \bigcup_{\ell \in
  \mathrm{supp}(p_{o_i})} \ell\right)\,.$$
\end{definition}

\begin{definition}
  Let $n\in\mathbbm{N}\setminus\{1,2\}$. We define an equivalence
  relation $\sim_{n}$ on $\mathbbm{R}^2$ by setting
$$
x\sim_{n} y\,\,\/\/ :\Longleftrightarrow\,\,\/\/ x-y \in \mathcal{O}_{n}\,.
$$
If $x,y\in\mathbbm{R}^2$ satisfy $x\sim_{n} y$, we say that $x$ and
$y$ are \emph{equivalent modulo} $\mathcal{O}_{n}$.
\end{definition}

\begin{definition}[Decomposition Problem]\label{decompset}
  Let $n\in\mathbbm{N}\setminus\{1,2\}$, and let $o_1,\dots,o_m$ be
  $m\geq 2$ pairwise non-parallel $\mathcal{O}_{n}$-directions. The
  corresponding decomposition problem is defined as follows.
  \begin{quote}
    {\sc Decomposition}. \\
    Given functions $p_{o_{i}} : \mathcal{L}_{o_{i}} \longrightarrow
    \mathbbm{N}_{0}$, $i\in\{1,\dots,m\}$, whose supports are finite
    and satisfy $\operatorname{supp}(p_{o_{i}})\subset
    \mathcal{L}^{\mathcal{O}_{n}}_{o_{i}}$, compute the equivalence
    classes modulo $\mathcal{O}_{n}$ in the associated grid
    $G_{\{p_{o_{i}}|i\in\{1,\dots,m\}\}}$.
\end{quote}
\end{definition}

Of course, this problem can be reduced to a polynomial number of
membership tests in $\mathcal{O}_{n}$. It is, however, not directly
clear how these tests can be performed and, actually, the algebraic
properties of $\mathcal{O}_{n}$ will be utilized. Also, later a
uniform bound for the number of classes will be given that is
independent of the $X$-ray data.

\begin{figure}
\centerline{\epsfysize=0.22\textheight\epsfbox{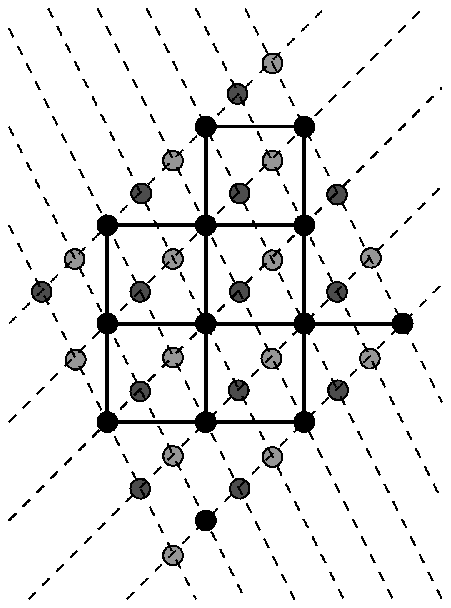}\hspace{0.2\textwidth}
\epsfysize=0.3\textheight\epsfbox{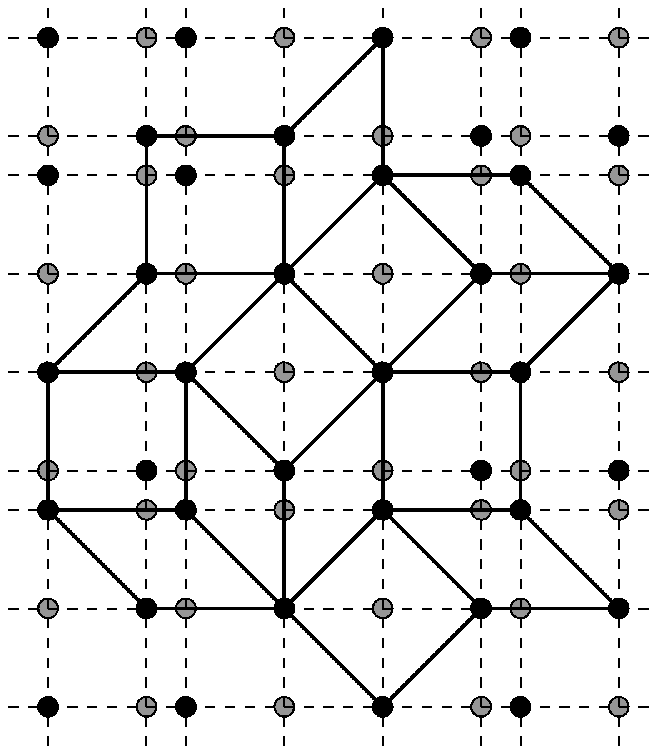}}
\caption{Grids arising from two $\mathcal{O}_{n}$-directions: On the left, the grid generated 
by the $X$-rays of a finite subset of a translate of $\mathcal{O}_{4}=\mathbbm{Z}^2$ 
in the two non-parallel $\mathcal{O}_{4}$-directions $(1,1)$ and $(1,-2)$. 
The three equivalence classes modulo $\mathcal{O}_{4}$ are marked by different greyscales. 
On the right, the grid generated by the $X$-rays of a finite subset of a translate 
of $\varLambda_{\rm AB}$ in the two non-parallel $\mathcal{O}_{8}$-directions $1$ and $\zeta_{8}^2=i$. 
The two equivalence classes are also shown.}
\label{fig:aequiv}
\end{figure}

\begin{remark}\label{m-rem:Zclasses}
  The phenomenon of multiple equivalence classes modulo
  $\mathcal{O}_{n}$ in the grid occurs already in the classical
  lattice situation; see Figure~\ref{fig:aequiv} on the left. There,
  \emph{no} translate of the finite subset of the square lattice
  (marked by the connecting lines) is contained in any of the other
  equivalence classes. Also, note the fact that \emph{exactly} one of
  the three equivalence classes has $14$ elements (the equivalence
  class marked by light grey), whereas the remaining two only have
  $13$ elements; it follows that this equivalence class (which
  generates the same grid as the marked finite subset of the square
  lattice) would be the unique solution of the corresponding
  reconstruction problem associated with its point set. Hence, the
  problem of decomposing the grid into its equivalence classes modulo
  $\mathcal{O}_{n}$ is the first problem to be solved when dealing
  with the consistency or the reconstruction problem, \emph{also} in
  the classical planar setting.
\end{remark}

\subsection{The Separation Problem}\label{m-subsec:sepprobl}

When dealing with the consistency, reconstruction and uniqueness
problems defined above, it is clear from the definition of
$W^{\circ}_{\mathcal{M}(\mathcal{O}_{n}),\star}$-sets that, given
$n\in\mathbbm{N}\setminus\{1,2,3,4,6\}$, a finite set $F$ of points in
$(\mathbbm{R}^2)^{\phi(n)/2-1}$ and a window
$W\subset(\mathbbm{R}^2)^{\phi(n)/2-1}$, we have to be able to decide
whether $F$ is contained in a translate of $W^{\circ}$. This leads us
to the following geometric separation problem for sets $[F]^{\star}
\subset t+W$.

\begin{definition}
  Let $d\in\mathbbm{N}$, let $P,W\subset\mathbbm{R}^d$, and let
  $t\in\mathbbm{R}^d$. We set
$$S_{W,t}(P):= P\cap (t+W)$$and, further,
 $$\operatorname{Sep}_{W}(P):= \left\{S_{W,t}(P)\,|\,t\in
  \mathbbm{R}^d\right\}\,.$$
\end{definition}

\begin{definition}[Separation Problem]\label{m-def:SWandSep}\hfill
  \begin{quote}
    {\sc Separation}. \\
    Given a finite set $P\subset \mathbbm{R}^d$, and a set $W\subset
    \mathbbm{R}^d$, determine $\operatorname{Sep}_W(P)$.
\end{quote}
\end{definition}

\begin{figure}
\begin{minipage}{\textwidth}
\begin{center}
\psfrag{W}{$W$}
\psfrag{-W}{$-\!W$}
\psfrag{t}{$t$}
\psfrag{p1}{$p_1$}
\psfrag{p2}{$p_2$} 
\psfrag{p3}{$p_3$}
 \psfrag{p1W}{$p_1\!\!-\!\!W$}
\psfrag{p2W}{$p_2\!\!-\!\!W$}
\psfrag{p3W}{$p_3\!\!-\!\!W$}
\includegraphics[scale=.55]{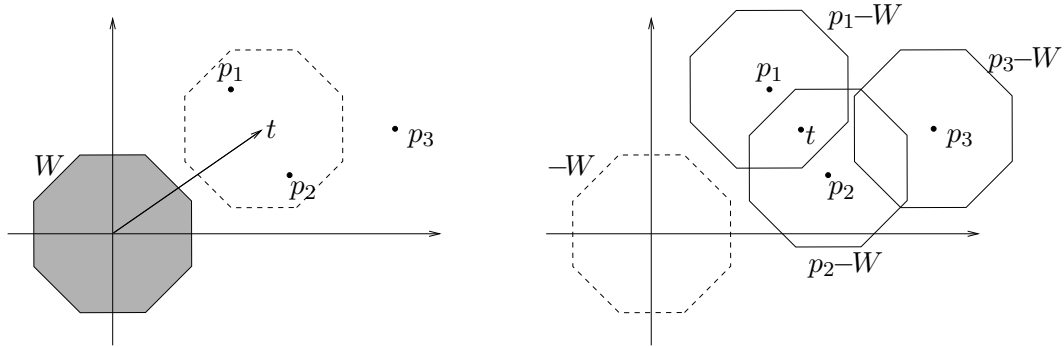}
\end{center}
\end{minipage}
\caption{On the left: If we translate $W$ by $t$, then $\{p_1,p_2\}$
  is a subset of $t+W$, but $\{p_3\}$ is not. On the right: The `world
  of translation vectors'. The point $t$ is contained in $p_1-W$ and
  $p_2-W$, but not in $p_3-W$. Again, we see that $S_{W,t}(P)=\{p_1,p_2\}$.
  (In this example, the window is centrally symmetric with respect to the origin, i.e., $W=-W$.)
}
\label{fig:changeview}
\end{figure}

\begin{remark}
  Note that $\operatorname{Sep}_{W}(P)$ contains all subsets of $P$
  that are `separable' from their complement (in $P$) by a translate
  of $W$. Trivially, one has $p\in t+W$ iff $t\in p-W$. It follows
  that
\begin{equation}\label{m-eq:ind}
S_{W,t}(P)=\{p\in P\,|\,t\in p-W\}\,.
\end{equation}
We will frequently make use of the above equivalence, because it
allows us to switch between a separable set $S_{W,t}(P)$ and the set
of translation vectors that makes it separable; see Fig.
\ref{fig:changeview} for an illustration.
\end{remark}

\subsection{Main Algorithmic Results}\label{mainresults}

In the following we apply the real RAM-model of computation, see
{\it e.g.}~Preparata \& Shamos (1985). Here each of the standard elementary
operations on reals counts only with unit cost.

Our first result shows that the decomposition problem can be solved
efficiently.

\begin{theorem}\label{th1}
  The decomposition problem can be solved in polynomial time in the
  real RAM model.  More precisely, it is of complexity $O(s^2)$, where
  $s$ is the maximum of the cardinalities of the supports of the given
  $X$-ray data functions.
\end{theorem}

The next result deals with the separation problem.

\begin{theorem}\label{separation}
  Let the window $W$ be given as an intersection of finitely many
  halfspaces, i.e., $W=\{x~|~Ax \leq b\}$ with
  $A\in\mathbbm{R}^{l\times d}$ and $b\in \mathbbm{R}^{l}$.  $($The
  parameters $d$, $A$, and $b$ are not part of the input$)$.  Then,
  for any finite set $P\subset\mathbbm{R}^d$, the problem of computing
  $\operatorname{Sep}_{W^{\circ}}(P)$ can be solved in
  $O(\operatorname{card}(P)^{d+1})$ operations.
\end{theorem}

As a consequence of Theorems \ref{th1} and \ref{separation} we see
that the standard tomographic algorithms that have been developed for
the lattice case can also be extended to the tomography of cyclotomic
model sets.

\begin{theorem}\label{easy}
  Let $W$ be given as in Theorem \ref{separation}. Then the problems
  {\sc Consistency}, {\sc Reconstruction} and {\sc Uniqueness} as
  defined in Definition \ref{m-def:consrecprob} can be solved with
  polynomially many operations and polynomially many calls to an
  oracle that solves the same problem on subsets of the plane of
  cardinality $O(s^2)$, where $s$ is again the maximum of the
  cardinalities of the supports of the given $X$-ray data functions.
\end{theorem}

As a simple corollary we finally note that the case of two directions
can be solved in polynomial time even for cyclotomic model sets.

\begin{corollary}\label{th2}
  When restricted to two $\mathcal{O}_{n}$-directions and polytopal
  windows the problems {\sc Consistency}, {\sc Reconstruction} and
  {\sc Uniqueness} as defined in Definition \ref{m-def:consrecprob}
  can be solved in polynomial time in the real RAM-model.
\end{corollary}

\section{Analysis of the Problems, Proofs and More Results} \label{m-sec:results}

In the following we give a detailed analysis of the problems
introduced in the previous section, prove the assertions stated there
and obtain more results on the way.

\subsection{Tractability of the Decomposition Problem}

We will now show that the number of equivalence classes of a grid is
uniformly bounded by a number that depends on the given directions but
is independent of the $X$-ray data. This result will then allow us to
prove Theorem \ref{th1}.

\begin{definition}\label{m-ref:compgrid}
  Let $n\in\mathbbm{N}\setminus\{1,2\}$ and let $o_1,o_2$ be two
  non-parallel $\mathcal{O}_{n}$-directions. We define the {\em
    complete grid} $G_{\{o_1,o_2\}}$ as
$$G_{\{o_1,o_2\}} := \bigcap_{i=1}^{2}\,\,\Big( \bigcup_{\ell \in \mathcal{L}^{\mathcal{O}_{n}}_{o_{i}}} \ell\Big)\,.$$    
\end{definition}

\begin{proposition}\label{index}
  Let $n\in\mathbbm{N}\setminus\{1,2\}$ and let $o_1,o_2$ be two
  non-parallel $\mathcal{O}_{n}$-directions. Then, the complete grid
  $G_{\{o_1,o_2\}}$ satisfies\/ $\mathcal{O}_{n}\subset
  G_{\{o_1,o_2\}}\subset \mathbbm{C}$ and $G_{\{o_1,o_2\}} \subset
  M_{\{o_1,o_2\}}$, where one sets
\begin{equation}\label{indexeq}
  M_{\{o_1,o_2\}}
  \ := \ 
  \operatorname{lin}_{\thinspace\scriptstyle{\mathcal{O}}\displaystyle_{n}}
  \left(\left\{\frac{1}{\alpha \delta - \beta \gamma} o_1,
      \frac{1}{\alpha \delta - \beta \gamma} o_2 \right\}\right)
  \,, 
\end{equation}
and $\alpha, \beta, \gamma , \delta \in
\thinspace\scriptstyle{\mathcal{O}}\displaystyle_{n}$ are determined
by $o_1 = \alpha + \beta \zeta_{n}$ and $o_2 = \gamma + \delta
\zeta_{n}$.
\end{proposition}

\begin{remark}\label{r3}
  Note that the linear independence of $\{o_1, o_2\}$ and $\{1,
  \zeta_{n}\}$ over $\mathbbm{R}$ implies that $\alpha \delta - \beta
  \gamma \neq 0$. Also, by definition $M_{\{o_1,o_2\}}$ is an
  $\thinspace\scriptstyle{\mathcal{O}}\displaystyle_{n}$-module of
  rank $2$ with basis\linebreak $\{o_1/(\alpha \delta - \beta \gamma),
  o_2/(\alpha \delta - \beta \gamma) \}$, and Proposition~\ref{index}
  shows that
$$\mathcal{O}_{n} \subset G_{\{o_1,o_2\}} \subset M_{\{o_1,o_2\}}
\subset \mathbbm{K}_{n}\,.$$ Note further that there are examples
where the inclusion $G_{\{o_1,o_2\}} \subset M_{\{o_1,o_2\}}$ is not
an equality. This is due to the fact that $M_{\{o_1,o_2\}}$ depends on
the scaling of $o_1$ and $o_2$, while $G_{\{o_1,o_2\}}$ does not. On
the other hand, let
$\gamma\in\thinspace\scriptstyle{\mathcal{O}}\displaystyle_{n}$ and
consider the two non-parallel $\mathcal{O}_{n}$-directions $1$ and
$\gamma+\zeta_{n}$. Then, by ~(\ref{indexeq}) and Lemma~\ref{Oo}, one
has
$$\mathcal{O}_{n}\subset G_{\{1,\gamma+\zeta_{n}\}} \subset
M_{\{1,\gamma+\zeta_{n}\}}=\mathcal{O}_{n}$$ and hence
$G_{\{1,\gamma+\zeta_{n}\}}=M_{\{1,\gamma+\zeta_{n}\}}=\mathcal{O}_{n}$.
Further observe that, for $n\in\mathbbm{N}\setminus\{1,2,3,4,6\}$, the
complete grid $G_{\{o_1,o_2\}}$ is a dense subset of the plane,
because already its subset $\mathcal{O}_{n}$ has this property; {\it cf.}
Remark~\ref{r0}.
\end{remark}

\begin{proof}[Proof of Proposition~\ref{index}]
The first inclusion is obvious by definition. 

Next, we claim that $\mathcal{O}_{n} \subset M_{\{o_1,o_2\}}$. Let $o
\in \mathcal{O}_{n}$.  By Lemma~\ref{Oo}(a), there are unique
$\varphi,\psi \in
\thinspace\scriptstyle{\mathcal{O}}\displaystyle_{n}$ with $o =
\varphi + \psi \zeta_{n}$.  By the linear independence of
$\{o_1,o_2\}$ over $\mathbbm{R}$, there are unique $x, y \in
\mathbbm{R}$ with $xo_1 + yo_2 = o$.  Hence
$$(x\alpha+y\gamma - \varphi) + (x\beta + y\delta - \psi) \zeta_{n} = 0$$ 
and, using the linear independence of $\{1,\zeta_{n}\}$ over $\mathbbm{R}$, 
we get that $x\alpha+y\gamma - \varphi = x\beta + y\delta - \psi = 0$. 
In matrix notation, this means that
$$\left(\begin{array}{cc}
\alpha & \gamma \\
\beta & \delta \\
\end{array}\right)
\left(\begin{array}{c}
x\\
y\\
\end{array}\right) = 
\left(\begin{array}{c}
\varphi\\
\psi\\
\end{array}\right)\,.$$  
Cramer's rule now implies that $$x = (\varphi \delta - \psi
\gamma)/(\alpha \delta - \beta \gamma) \in
\thinspace\scriptstyle{\mathcal{O}}\displaystyle_{n}/(\alpha \delta -
\beta \gamma)$$ and $$y = (\alpha \psi - \beta \varphi)/(\alpha \delta -
\beta \gamma) \in
\thinspace\scriptstyle{\mathcal{O}}\displaystyle_{n}/(\alpha \delta -
\beta \gamma)\,.$$  This proves our claim.

Finally, consider $g \in G_{\{o_1,o_2\}}$. By definition, there are
elements 
$o' , o'' \in \mathcal{O}_{n}$ such that $\{g\} = (o' + \mathbbm{R} o_1)
\cap (o'' + \mathbbm{R} o_2)$. Moreover, there are unique $x,y \in
\mathbbm{R}$ with $g = o' + x o_1 = o'' + y o_2$.  Hence, $x o_1 +
(-y) o_2 = o'' - o' \in \mathcal{O}_{n}$ and, by the same calculation
as above, we get that $x, y \in
\thinspace\scriptstyle{\mathcal{O}}\displaystyle_{n}/(\alpha \delta -
\beta \gamma)$.  Together with our first claim, this shows that $g \in
M_{\{o_1,o_2\}}$.
\end{proof}

\begin{lemma}\label{phin}
$M_{\{o_1,o_2\}}$ is a $\mathbbm{Z}$-module of rank $\phi(n)$.
\end{lemma}
\begin{proof}
  This is an immediate consequence of the facts that $M_{\{o_1,o_2\}}$
  is an $\thinspace\scriptstyle{\mathcal{O}}\displaystyle_{n}$-module
  of rank~$2$ and
  $\thinspace\scriptstyle{\mathcal{O}}\displaystyle_{n}$ is a
  $\mathbbm{Z}$-module of rank $\phi(n)/2$; see Remark~\ref{r3} and
  Remark~\ref{r7}.
\end{proof} 

The following lemma shows that $M_{\{o_1,o_2\}}$, and thus
$G_{\{o_1,o_2\}}$, decomposes into finitely many equivalence classes
whose number depends only on $\{o_1,o_2\}$. Note that the symbol
$\dot{\cup}$ is used to indicate disjoint unions.

\begin{lemma}\label{MOindex}
  The subgroup index $[M_{\{o_1,o_2\}} : \mathcal{O}_{n}]$ is finite.
  Hence, there are $c \in \mathbbm{N}$ and $t_1, t_2, \dots, t_c $
  $\in M_{\{o_1,o_2\}}$ such that $$M_{\{o_1,o_2\}} =
  \dot{\bigcup}_{i=1}^{c} (t_{i} + \mathcal{O}_{n})\,,$$ where,
  without restriction, $t_1=0$. It follows that every subset $G$ of
  $M_{\{o_1,o_2\}}$ satisfies the decomposition
$$G = \dot{\bigcup}_{i=1}^{c} (G \cap (t_{i} + \mathcal{O}_{n}))\,.$$ 
\end{lemma}
\begin{proof}
  By Lemma~\ref{phin}, $M_{\{o_1,o_2\}}$ is a $\mathbbm{Z}$-module of
  rank $\phi(n)$.  Moreover, $M_{\{o_1,o_2\}}$ is torsion-free because
  it is a subset of the field $\mathbbm{K}_{n}$; see Remark~\ref{r3}.
  But $\mathcal{O}_{n}$ is a $\mathbbm{Z}$-module of rank $\phi(n)$ as
  well; see Remark~\ref{r7}. Now, Proposition~\ref{full} yields the
  result.
\end{proof}

\begin{remark}
  By Proposition~\ref{full}, the subgroup index $$[M_{\{o_1,o_2\}} :
  \mathcal{O}_{n}]$$ equals the absolute value of the determinant of
  the transition matrix $A$ from any\/ $\mathbbm{Z}$-basis of
  $M_{\{o_1,o_2\}}$ to any\/ $\mathbbm{Z}$-basis of $\mathcal{O}_{n}$.
  It follows that, given the $\mathbbm{Z}$-coordinates of $o_1$ and
  $o_2$ with respect to the $\mathbbm{Z}$-basis
  $\{1,\zeta_{n},\zeta_{n}^2, \dots ,\zeta_{n}^{\phi(n)-1}\}$ of
  $\mathcal{O}_{n}$ ({\it cf.} Remark~\ref{r7}), one is able to compute
  $[M_{\{o_1,o_2\}} : \mathcal{O}_{n}]$. Note that, for any
  $\gamma\in\thinspace\scriptstyle{\mathcal{O}}\displaystyle_{n}$, one
  has $$[M_{\{1,\gamma+\zeta_{n}\}} : \mathcal{O}_{n}]=1\,;$$see
  Remark~\ref{r3}.
\end{remark}

\begin{remark}\label{r4}  Let $n\in\mathbbm{N}\setminus\{1,2\}$ and let
  $o_1,\dots,o_m$ be $m\geq 2$ pairwise non-parallel
  $\mathcal{O}_{n}$-directions. For any instance of the corresponding
  decomposition problem, the associated grid
  $G_{\{p_{o_{i}}|i\in\{1,\dots,m\}\}}$ satisfies
$$\operatorname{card}(G_{\{p_{o_{i}}|i\in\{1,\dots,m\}\}}) \leq s^2\,,$$ where
$$s:=\operatorname{max}(\{ \operatorname{card}(\operatorname{supp}(p_{o_{i}}))
\, |\, i\in\{1,\dots,m\}\})\,.$$ Since
$G_{\{p_{o_{i}}|i\in\{1,\dots,m\}\}} \subset G_{\{o_1,o_2\}}$,
Proposition~\ref{index} shows that the last part of
Lemma~\ref{MOindex} applies to $G_{\{p_{o_{i}}|i\in\{1,\dots,m\}\}}$.
\end{remark}
 
In the following we assume that the elements of the supports of the
$p_{o_{i}}$, $i\in\{1,\dots,m\}$, are given in the form $o +
\mathbbm{R} o_{i}$ for suitable $o \in \mathcal{O}_{n}$. Moreover, we
assume that all $o$'s and all $o_{i}$ are given through their
$\mathbbm{Z}$-coordinates with respect to the $\mathbbm{Z}$-basis
$\{1,\zeta_{n},\zeta_{n}^2, \dots ,\zeta_{n}^{\phi(n)-1}\}$ of
$\mathcal{O}_{n}$ ({\it cf.} Remark~\ref{r7}).

We now prove Theorem \ref{th1} which we restate in a rephrased
form.

\begin{theorem}\label{th1a}
  The decomposition problem can be solved with $O(s^2)$ many real
  number operations.
\end{theorem}

\begin{proof} The algorithm performs the following steps.

  \textit{Step 1:} By the proof of Lemma~\ref{Oo}(a), the Euclidean
  algorithm in $\mathbbm{Z}[X]$, the inductive computability of the
  $n$th cyclotomic polynomial
  $F_{n}=\operatorname{Mipo}_{\mathbbm{Q}}(\zeta_{n})$ ({\it cf.}
  Remark~\ref{recur} and Proposition~\ref{gau2}), the proof of
  Proposition~\ref{index} and the Gaussian elimination algorithm, we
  are able to compute the $\mathbbm{Q}$-coordinates of the elements of
  the grid $G_{\{p_{o_{i}}|i\in\{1,\dots,m\}\}}\subset
  \mathbbm{K}_{n}$ with respect to the $\mathbbm{Q}$-basis
$$\{1,\zeta_{n},\zeta_{n}^2, \dots ,\zeta_{n}^{\phi(n)-1}\}$$ 
of $\mathbbm{K}_{n}$ ({\it cf.} Proposition~\ref{gau}) efficiently.\medskip

\textit{Step 2:} Since $\{1,\zeta_{n},\zeta_{n}^2, \dots
,\zeta_{n}^{\phi(n)-1}\}$ is simultaneously a $\mathbbm{Q}$-basis of
$\mathbbm{K}_{n}$ and a $\mathbbm{Z}$-basis of $\mathcal{O}_{n}$ ({\it cf.}
Proposition~\ref{gau} and Remark~\ref{r7}), one has for all $q_{0},
q_{1}, \dots, q_{\phi(n) - 1} \in \mathbbm{Q}$ the equivalence
\begin{equation}\label{criterion}
\begin{array}{lcl}
&&q_{0} + q_{1} \zeta_{n} + \dots + q_{\phi(n) - 1} \zeta_{n}^{\phi(n) - 1} 
\in \mathcal{O}_{n}\vspace{1ex}\\ &\Longleftrightarrow&   q_{0}, q_{1}, \dots,
q_{\phi(n) - 1} \in \mathbbm{Z}\,.
\end{array}
\end{equation}
By Step 1, the elements of $G_{\{p_{o_{i}}|i\in\{1,\dots,m\}\}}$ are
given in the form
$$q_{0} + q_{1} \zeta_{n} + \dots + q_{\phi(n) - 1} \zeta_{n}^{\phi(n) - 1}\,,$$ 
where $q_{0}, q_{1}, \dots, q_{\phi(n) - 1} \in \mathbbm{Q}$.  Now,
proceed as follows: choose an arbitrary element $g$ of
$G_{\{p_{o_{i}}|i\in\{1,\dots,m\}\}}$ and compute the
$\mathbbm{Q}$-coordinates of the differences $g - h$ with respect to
$\{1,\zeta_{n},\zeta_{n}^2, \dots ,\zeta_{n}^{\phi(n)-1}\}$, where $h
\in G_{\{p_{o_{i}}|i\in\{1,\dots,m\}\}} \setminus \{g\}$.  By the
above criterion~(\ref{criterion}), a fixed $h$ lies in the same
equivalence class modulo $\mathcal{O}_{n}$ as $g$ iff all coordinates
of $g - h$ are elements of $\mathbbm{Z}$. Iterate this procedure by
successively removing the computed equivalence classes and proceeding
with the remaining subset of the grid and an arbitrary element
therein.

We already saw in Remark~\ref{r4} that the last part of
Lemma~\ref{MOindex} applies to $G_{\{p_{o_{i}}|i\in\{1,\dots,m\}\}}$.
This immediately implies that Step 2 of this algorithm computes the
equivalence classes of the grid modulo $\mathcal{O}_{n}$ in at most
$$c:=[M_{\{o_1,o_2\}} : \mathcal{O}_{n}]\in\mathbbm{N}$$ 
iterations. The inequality $\operatorname{card}(G_{\{p_{o_{i}}|i\in\{1,\dots,m\}\}}) \leq s^2$ 
({\it cf.} Remark~\ref{r4}) now completes the proof.
\end{proof}

\begin{remark}
  The proof of Theorem \ref{th1a} indicates that we actually do not
  need the full strength of the real RAM-model of computation. Rather,
  a Turing machine model that is augmented for algebraic computations
  suffices, see {\it e.g.} Buchberger {\it et al.} (1982). Then, of course, the
  binary size of the input matters.
\end{remark}

\subsection{Tractability of the Separation Problem}
    
The problem {\sc Separation} in its general form is interesting on its
own and we show now how to deal with it for windows $W$ that are open
polyhedra, i.e.,
\[W=\{x~|~Ax < b\} \quad \text{with $A\in\mathbbm{R}^{l\times d}$ and
  $b\in \mathbbm{R}^{l}$}\;,\] where $d\geq 2$ is a fixed constant.
The ideas presented here can be generalized to semialgebraic sets, but
we prefer to keep the exposition more elementary.  Also, polytopal
windows with $N$-fold cyclic symmetry, where $N$ is the function from
(\ref{eq}), are most relevant for model sets. (Note that the windows
underlying the examples in Section \ref{m-sec:examples} are
polytopes.)

We will begin with some standard facts about hyperplane arrangements
as they are needed to deal with {\sc Separation}.  See Edelsbrunner
{\it et al.} (1986) for more information on hyperplane arrangements, and
 Agarwal \& Sharir (2000) and Halperin (2004) for surveys that cover
also more general classes of arrangements.

\begin{definition}
  For $i\in\{1,\dots,l\}$, let $a_i\in\mathbbm{R}^d\setminus\{0\}$,
  $\beta_i\in\mathbbm{R}$, and consider the sets
  $H_i=\{x~|~a_i^Tx=\beta_i\}$. Then $H_i$ is called {\em hyperplane}
  and $\mathcal{H}=\{H_1,\dots,H_{l}\}$ is a \emph{hyperplane
    arrangement} in $\mathbbm{R}^d$.  The \emph{sign vector} $SV(x)$
  of some point $x\in\mathbbm{R}^d$ is defined component-wise via
\[SV_i(x):=\begin{cases}
-1 & \text{if $a_i^Tx<\beta_i$}\\
~0 & \text{if $a_i^Tx=\beta_i$}\\
+1 & \text{if $a_i^Tx>\beta_i$}
\end{cases}~~\quad,~~ 1\leq i\leq l~~.\]
For $s\in\{\pm 1,0\}^l$ with
\[C_{s}:=\{x~|~SV(x)=s\}\ne \varnothing \, ,\]
$C_{s}$ is called a $($proper$)$ \emph{cell} of the arrangement $\mathcal{H}$.
\end{definition}

\begin{figure}
\begin{minipage}{\textwidth}
\begin{center}
\psfrag{+}{$+$}
\psfrag{-}{$-$}
\psfrag{h1}{$H_1$}
\psfrag{h2}{$H_2$}
\psfrag{h3}{$H_3$}
\psfrag{h4}{$H_4$}
\psfrag{s1}{\scalebox{0.8}{$(0,0,+1,-1)$}}
\psfrag{s2}{\scalebox{0.8}{$(+1,0,+1,-1)$}}
\psfrag{s3}{\scalebox{0.8}{$(-1,-1,-1,+1)$}}
\includegraphics[scale=.55]{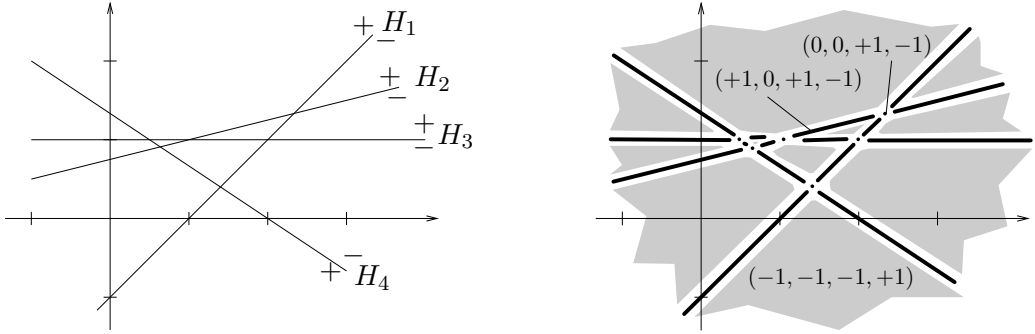}
\end{center}
\end{minipage}
\caption{A small example for a hyperplane-arrangement in the plane.
The hyperplanes are given by $H_i=\{x~|~a_i^Tx=\beta_i\}$, where
$a_1=(-1,1)$, $\beta_1=-1$,
$a_2=(-1,4)$, $\beta_2=3$,
$a_3=(0,1)$, $\beta_3=1$,
$a_4=(-2,-3)$, $\beta_4=-4$.
On the right, the cells are drawn schematically.
The arrangement consists of six points, 16 one-dimensional cells (thick lines) 
and 11 full-dimensional cells (grey areas). Some sign vectors are given. 
Note that not all vectors in $\{\pm 1,0\}^4$ occur as sign vectors of cells; 
{\it e.g.}, $(0,0,0,0)$ and $(-1,+1,-1,-1)$ are not realized.
}
\label{fig:arrexa}
\end{figure}

\begin{remark}
  The cells of an arrangement are relatively open sets of various
  dimensions.  In particular, a cell $C_{s}$ with sign vector $s$ is
  full-dimensional if and only if $s\in\{\pm 1\}^l$.  Of course,
  $\mathbbm{R}^d$ is the disjoint union of all the cells of a
  hyperplane arrangement.  Figure \ref{fig:arrexa} gives some
  illustration.
\end{remark}

In view of their general relevance, hyperplane arrangements are well
studied and also algorithmically well understood.

\begin{proposition}\label{m-prop:linear} 
  Let $\mathcal{H}=\{H_1,\dots,H_{l}\}$ be a hyperplane arrangement in
  $\mathbbm{R}^d$.  There exists an algorithm that computes a set of
  points meeting each cell of $\mathcal{H}$ in ${{O}}(l^{d})$
  operations in the real RAM model.
\end{proposition}

\begin{proof}
  {\it Cf.} Theorem 3.3 of Edelsbrunner {\it et al.} (1986). See also Chapter 7
  of Edelsbrunner (1987).\end{proof}

The proof of Theorem \ref{separation} will now be based on the
following observation that ties the separation problem to certain
hyperplane arrangements.

\begin{figure}
\begin{minipage}{\textwidth}
\begin{center}
\psfrag{W}{$W$}
\psfrag{-W}{$-\!W$}
\psfrag{t}{$t$}
\psfrag{t'}{$t'$}
\psfrag{p1}{$p_1$}
\psfrag{p2}{$p_2$}
\psfrag{p3}{$p_3$}
 \psfrag{p1W}{$p_1\!\!-\!\!W$}
\psfrag{p2W}{$p_2\!\!-\!\!W$}
\psfrag{p3W}{$p_3\!\!-\!\!W$}
\includegraphics[scale=.55]{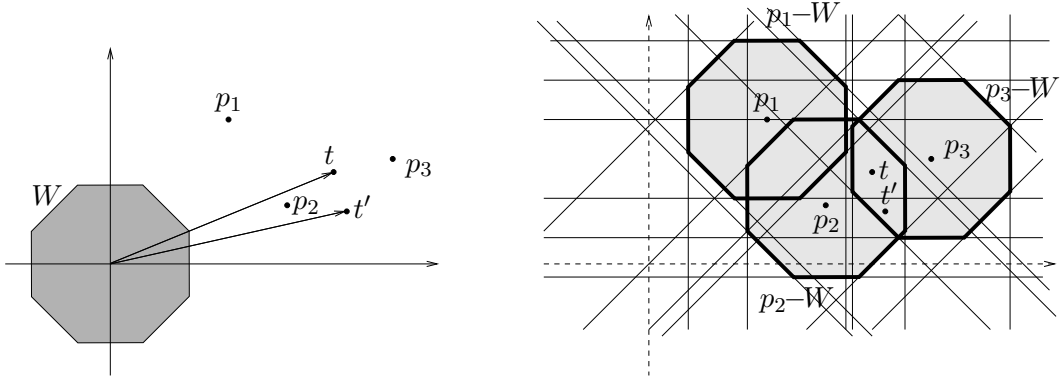}
\end{center}
\end{minipage}
\caption{This shows in essence the same situation as in Fig. \ref{fig:changeview}, on the right we added the arrangement $\mathcal{H}(W,P)$. Here, $S_{W,t}(P)=S_{W,t'}(P)$, but $t$ and $t'$ are in different fulldimensional cells of the arrangement $\mathcal{H}(W,P)$, see Lemma \ref{m-lem:fullarr}.
Therefore the inverse direction of the implication in Lemma \ref{m-lem:fullarr} (b) is not true.}
\label{fig:samesep}
\end{figure}

\begin{lemma} \label{m-lem:fullarr} Let $P=\{p_1,\dots,p_q\}$ be a
  finite set of points in $\mathbbm{R}^d$, let $W=\{x~|~Ax< b\}$ with
  $A\in\mathbbm{R}^{l\times d}$, $b\in \mathbbm{R}^{l}$, and let
  $a_i^T$ denote the $i$th row of $A$, $1\leq i\leq l$. For $1\leq
  i\leq l$ and $1\leq j\leq q$, set
$$
H^{(j)}_{i}:=\{x~|~a_i^Tx=(Ap_j-b)_i\} 
$$
Further, set
$$
\mathcal{H}(W,P):=\{H^{(j)}_{i}~|~1\leq i\leq l ~,~1\leq j\leq q\}\;.
$$
Then, one has the following:
\begin{enumerate}
\item[{\rm (a)}] The set $p_j-W$ is an intersection of open halfspaces
  defined by the $H^{(j)}_{1},\dots, H^{(j)}_{l}$, more precisely
\[p_j-W=\{x~|A^Tx > Ap_j-b\}\;.\]
\item[{\rm (b)}] For each cell $C_{s}$ of the hyperplane arrangement
  $\mathcal{H}(W,P)$ with sign vector $s=(s_{i,j})_{i,j}$ the
  following implication is true:
  \[t,t'\in C_{s}\quad \Longrightarrow\quad
  S_{W,t}(P)=S_{W,t'}(P)\;.\] $($Of course, the reverse implication is
  not true in general; see Figure \ref{fig:samesep}.$)$
\end{enumerate}
\end{lemma}

\begin{proof}
Part (a) follows from a simple computation:
\begin{eqnarray*} 
p-W\,=\,\{p-x~|~Ax< b\}&=&\{x~|~A(p-x)< b\}\\
&=&\{x~|~Ax>Ap-b\}\;.
\end{eqnarray*}
For (b), recall from (\ref{m-eq:ind}) that
$$S_{W,t}(P)=\{p_j~|~1\leq j\leq q~,~t\in p_j-W\}$$ for any $t\in\mathbbm{R}^d$. 
Using (a) we see that $t\in p_j-W$ iff
$SV_{1j}(t)=\dots=SV_{lj}(t)=+1$. Now, if $t,t'\in C_{s}$, we have
$SV(t)=SV(t')=s$, concluding the proof.
\end{proof}

Here is a restatement of Theorem \ref{separation}.

\begin{theorem}\label{m-thm:seplinear}
  Let $W=\{x~|~Ax< b\}$ with $A\in\mathbbm{R}^{l\times d}$ and $b\in
  \mathbbm{R}^{l}$. Moreover, let $P=\{p_1,\dots,p_q\}$ be a finite
  set of points in $\mathbbm{R}^d$.  Then, $\operatorname{Sep}_W(P)$
  can be computed with the aid of at most ${{O}}((lq)^{d+1})$
  operations in the real RAM model.
\end{theorem}

\begin{proof}
  Our algorithm to determine $\operatorname{Sep}_W(P)$ performs the following steps.\\
  \textit{Step 1:} Compute $(Ap_j-b)_i$ for $1\leq i\leq l$ and $1\leq
  j\leq q$, to
  specify the hyperplane arrangement $\mathcal{H}(W,P)$ from Lemma \ref{m-lem:fullarr}.\\
  \textit{Step 2:} Compute a set $T$ of points meeting every cell of $\mathcal{H}(W,P)$.\\
  \textit{Step 3:} For each of the points $t\in T$ obtained in 2, compute $S_{W,t}(P)$.\\
  \textit{Step 4:} Output the collection of all the $S_{W,t}(P)$.

  The correctness of this procedure follows directly from Lemma
  \ref{m-lem:fullarr}.

  Now we show the complexity assertion.  Step 1 needs no more than
  ${{O}}(lq)$ operations.  Step~2 requires ${{O}}((lq)^d)$ operations
  by Proposition \ref{m-prop:linear}.  For Step 3, we decide if $t\in
  p_j-W$ for each $j$. To this end we test if $t$ satisfies the
  inequalities $a_i^Tt > (Ap_j-b)_i$, $1\leq i\leq l$, $1\leq j\leq
  q$.  This is done with ${{O}}(lq)$ operations.  In total we do not
  need more than ${{O}}(lq+(lq)^dlq)=O((lq)^{d+1})$ operations.
\end{proof}

\begin{remark}\label{rem-size} 
  As the proof of Theorem \ref{m-thm:seplinear} shows, if the number
  of hyperplanes defining the window $W$ is regarded constant, then
$$
\operatorname{card}\left(\operatorname{Sep}_{W^{\circ}}(P)\right)=O\left(\operatorname{card}(P)^{d}\right)\,.
$$
\end{remark}

\begin{remark} Theorem \ref{m-thm:seplinear} can be generalized to
  semialgebraic sets $W$. The corresponding algorithm is then based on
  an analogue of Proposition \ref{m-prop:linear} in the semialgebraic
  world; see Basu {\it et al.} (1996) and Theorem 2 of Basu {\it et al.} (1997).
\end{remark}

\subsection{On the Tractability of Consistency, Reconstruction and Uniqueness}

As a consequence of Theorems \ref{th1} and \ref{separation} we can now
prove Theorem \ref{easy}.  In the following we only deal with {\sc
  Consistency} in detail; the proofs for the other two problems are
similar.  As Theorem~\ref{easy} states, we want to reduce {\sc
  Consistency} to a problem in the classical (anchored) case.

The number $m\in \mathbbm{N}\setminus \{1\}$ of $X$-rays and the
different directions $o_{1},\dots,o_{m}$ are of course fixed as usual.

\begin{quote}
  {\sc AnchoredConsistency}. \\
  Given $s\in\mathbbm{N}$ and $p_{o_{i}} : \mathcal{L}_{o_{i}}
  \longrightarrow \mathbbm{N}_{0}$, $i\in\{1,\dots,m\}$, with finite
  supports whose cardinalites are bounded by $s$, and a finite set
  $S\subset \mathbbm{R}^2$ with at most $s^2$ points. Decide whether
  there is a set $F$ contained in $S$ which satisfies
  $X_{o_{i}}F=p_{o_{i}}$, $i\in\{1,\dots,m\}$.
\end{quote}

Now we show that for polytopal windows the problem {\sc Consistency}
for cyclotomic model sets can be reduced to {\sc AnchoredConsistency}.
Let $\mathcal{A}$ be an algorithm for solving {\sc
  AnchoredConsistency}. (In the following $\mathcal{A}$ acts as a
black box subroutine for the reduction.)

\begin{theorem}\label{easyconsistency}
  Let $W$ be given as in Theorem \ref{separation}. Then {\sc
    Consistency} can be solved with polynomially many operations and
  polynomially many calls to $\mathcal{A}$.
\end{theorem}

\begin{proof} The algorithm performs the following steps.

  \textit{Step 1:} Check first the necessary condition that the
  cardinalities
$$
\sum_{l \in
\operatorname{supp}(p_{o_{i}})}p_{o_{i}}(l)
$$
coincide for each $i$.
If this is the case, proceed with Step~2.
Otherwise the instance is inconsistent.
\medskip

\textit{Step 2:} Compute the elements of the equivalence classes
$G_{i}$ of the associated grid\linebreak $G_{\{ p_{o_{1}},\dots, p_{o_{m}} \}}$
modulo $\mathcal{O}_{n}$, say
$$G_{\{ p_{o_{1}}, \dots, p_{o_{m}} \}}=\dot\bigcup_{i=1}^{c}G_{i}\subset\mathbbm{K}_{n}$$
in terms of their $\mathbbm{Q}$-coordinates with respect to the
$\mathbbm{Q}$-basis
$$\{1,\zeta_{n},\zeta_{n}^2, \dots ,\zeta_{n}^{\phi(n)-1}\}$$ 
of $\mathbbm{K}_{n}$ ({\it cf.} Proposition~\ref{gau}). By
Theorem~\ref{th1}, this can be done efficiently. 
\medskip

\textit{Step 3:} For all $i\in\{1,\dots,c\}$, compute the
$.^{\star}$-image $[G_{i}]^{\star}$ of $G_{i}$. Note that we consider
the star map here as a map
 $$.^{\star}\! : \, \mathbbm{K}_{n}\longrightarrow (\mathbbm{R}^2)^{\frac{\phi(n)}{2}-1}\,.$$ 
 This can be done efficiently. Due to the definition of
 $W^{\circ}_{\mathcal{M}(\mathcal{O}_{n}),\star}$-sets, a solution
 $F\subset G_{i}$ for our instance must satisfy the condition
\begin{equation}\label{sepcond}
\exists~\tau\in(\mathbbm{R}^2)^{\frac{\phi(n)}{2}-1}~:~ [F]^{\star}\subset \tau+W^{\circ}\, . 
\end{equation}
Recall that for $n\in\{3,4,6\}$, condition~(\ref{sepcond}) is always
satisfied and one can proceed with Step~$4$. Otherwise, compute the
set $ \operatorname{Sep}_{W^{\circ}}([G_{i}]^{\star})$. By
Theorem~\ref{separation}, this can be done efficiently.  Note that,
for every $i\in\{1,\dots,c\}$, a subset $F\subset G_{i}$ that
satisfies condition~(\ref{sepcond}) has the property that
$[F]^{\star}\subset P$ for a suitable
$P\in\operatorname{Sep}_{W^{\circ}}([G_{i}]^{\star})$. Finally,
compute, for all $i\in\{1,\dots,c\}$ and for all
$P\in\operatorname{Sep}_{W^{\circ}}([G_{i}]^{\star})$, the pre-images
$S:=[P]^{-\star}$ of $P$ under the star map. This can be done
efficiently.  Note that, with the above restriction
$n\notin\{3,4,6\}$, the star map is injective. \medskip

\textit{Step 4:} If $n\in\{3,4,6\}$, consider the equivalence classes
$S:=G_{i}$, $i\in\{1,\dots,c\}$, having the property that
$\operatorname{card}(G_{i})\geq N$.  Otherwise, consider, for all
$i\in\{1,\dots,c\}$ and for all
$P\in\operatorname{Sep}_{W^{\circ}}([G_{i}]^{\star})$, the subsets
$S:=[P]^{-\star}$ of $G_{i}$ having the property that
$\operatorname{card}([S]^{-\star})\geq N$.  Then apply $\mathcal{A}$
on each such $S$. The instance is consistent iff $\mathcal{A}$ reports
consistency for one of the sets $S$.
\end{proof}

Note that for $m=2$ a polynomial-time algorithm $\mathcal{A}$ is
available; see {\it e.g.}~Slump \& Gerbrands (1982). There it is shown how
to set up a capacitated network that admits a certain flow iff the
consistency question has an affirmative answer. Points in the grid
correspond to arcs in this network. If we want to forbid certain
positions, we only have to cancel the corresponding arcs.  Hence we
obtain Corollary \ref{th2} for {\sc Consistency}.

The proofs for {\sc Reconstruction} and {\sc Uniqueness} are analogous.

\begin{figure}
\begin{minipage}{\textwidth}
\begin{center}
\psfrag{phy}{{\scriptsize physical space}}
\psfrag{int}{{\scriptsize internal space}}
\psfrag{star}{{\scriptsize$\,\,\,\,\,\,\,.^{\star}$}}
\psfrag{invstar}{{\scriptsize\parbox{1cm}{$\,\,\,\,\,\,\,.^{-\star}$}}}
\psfrag{1-1}{{\scriptsize 1--1}}
\psfrag{1}{$1$}
\psfrag{W}{$W$}
\includegraphics[scale=.55,clip]{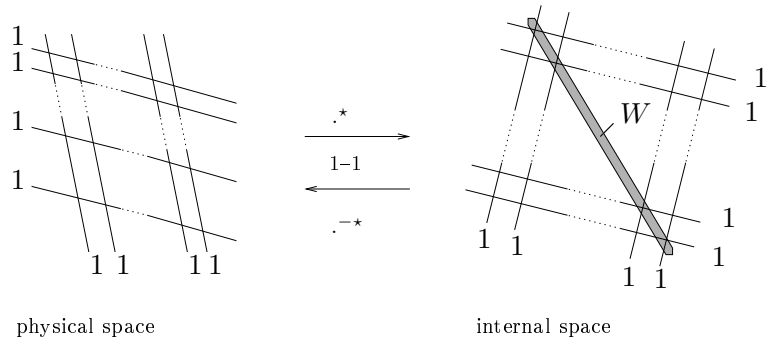}
\end{center}
\end{minipage}
\caption{A $(\operatorname{card}(F)$$\times$$\operatorname{card}(F)$)-grid (left hand side).
The subsets of the grid that conform to the $X$-ray data correspond to permutation matrices, hence there are $(\operatorname{card}(F))!$ of them. 
Assume that the internal space is also two-dimensional and that the
star map $.^{\star}$ acts as a bijection on the grid points. Then the grid in the physical space is mapped to a grid in the internal space (right hand side). 
(As concrete example of this situation, we recall the examples in
Section \ref{m-sec:examples}. {\it E.g.}, if $n=8$ and the star map is
defined by $\xi_8\mapsto \xi_8^3$, then the grid in the physical space
can be chosen as a finite patch of  
the $\mathbbm{Z}$-span 
$\operatorname{lin}_{\mathbbm{Z}}(\{1,\xi_8\})$
 (a lattice in $\mathbbm{R}^2$). It is in 1--1 correspondence to its
 image under the star map.
If the window is chosen `slim' enough, covering a `diagonal' of the
image grid, then there is only one of the $(\operatorname{card}(F))!${\ } solutions (in the internal space) that can be covered by a translate of the window. Thus, there is also only one single solution in the physical space.
}
\label{fig:firstsep}
\end{figure}

\begin{remark}\label{m-rem:firstsep}
  Note that the seemingly more natural approach to find subsets
  $F\subset G_{i}$ \emph{first} that conform to the $X$-rays, and
  check \emph{then} whether (\ref{sepcond}) is satisfied may lead to
  an exponential running time. In fact, Figure \ref{fig:firstsep}
  gives a simple example with a unique solution but exponentially many
  subsets of the grid conforming to the $X$-ray data.
\end{remark}

\section*{Acknowledgements}
It is a pleasure to thank U. Grimm for fruitful discussions. This work was supported by the German Research Council (Deutsche
Forschungsgemeinschaft), MB and CH within the CRC
701, and PG and KL through GRK 447/2.

\end{document}